\documentclass[11pt,a4paper,leqno]{article}
\usepackage{amsmath}
\usepackage{amssymb}
\usepackage{graphicx}
\advance\topmargin by -2.2cm
\newcommand{\bbr}{I\!\!R}

\newcommand{\bbn}{I\!\!N}
\newcommand{\bbz}{Z\!\!\!Z}

\newcommand{\calb}{{\cal B}}

\newcommand{\cald}{{\cal D}}
\newcommand{\cale}{{\cal E}}

\newcommand{\call}{{\cal L}}

\newcommand{\XXX}{{\bar X}}

\newcommand{\barr}{\begin{array}}
\newcommand{\earr}{\end{array}}
\newcommand{\beqq}{\begin{equation}}
\newcommand{\eeqq}{\end{equation}}
\newcommand{\beao}{\begin{eqnarray*}}
\newcommand{\eeao}{\end{eqnarray*}\noindent}
\newcommand{\beam}{\begin{eqnarray}}
\newcommand{\eeam}{\end{eqnarray}\noindent}

\newcommand{\halmos}{\quad\hfill\mbox{$\Box$}}

\newcommand{\la}{\lambda}

\newcommand{\si}{\sigma}
\newcommand{\al}{\alpha}
\newcommand{\vth}{\vartheta}

\newcommand{\vep}{\varepsilon}

\newtheorem{theo}{Theorem}

\newtheorem{cor}{\indent Corollary}

\newtheorem{ex}{\indent Example}
\newtheorem{ass}{\indent Assumption}

\newcommand{\wt}{\widetilde}

\newcommand{\ov}{\overline}

\newcommand{\lra}{\longrightarrow}

\author{R. H\"opfner$^*$ \and E.~L\"ocherbach \and M. Thieullen\thanks{This work has been supported by the Agence Nationale de la Recherche through the project MANDy, Mathematical Analysis of Neuronal Dynamics, ANR-09-BLAN-0008-01. e-mail addresses: \tt{hoepfner@mathematik.uni-mainz.de},
     \tt{eva.loecherbach@u-cergy.fr} and
    \tt{michele.thieullen@upmc.fr}} \\
{\it Johannes Gutenberg-Universit\"at Mainz, Universit\'e de Cergy-Pontoise} \\
 {\it and Universit\'e Pierre et Marie Curie.}}
\begin{document}

\title{Transition densities for stochastic Hodgkin-Huxley models}

\maketitle

\begin{abstract}
We consider a stochastic Hodgkin-Huxley model driven by a periodic signal as model for the membrane potential of a pyramidal neuron. The  associated five dimensional diffusion process is a time inhomogeneous highly degenerate diffusion for which the weak H\"ormander condition holds only locally. Using a technique which is based on estimates of the Fourier transform, inspired by Fournier 2008, Bally 2007 and De Marco 2011, we show that the process admits locally a strictly positive continuous transition density. Moreover, we show that the presence of noise enables the stochastic system to imitate any possible deterministic spiking behavior, i.e. mixtures of regularly spiking and non-spiking time periods  are possible features of the stochastic model. This is a fundamental difference between stochastic and deterministic Hodgkin-Huxley models.  
\end{abstract}

\\


\section{Introduction}
In this paper we study stochastic models based on the Hodgkin-Huxley model well-known in physiology. Our focus is on the presence of a periodic stochastic stimulus standing for the synaptic input received by a neuron from a large number of other neurons. This leads us to the study of a highly degenerate non time homogeneous stochastic system which can not be addressed by available techniques from the literature. 
  
The deterministic Hodgkin-Huxley model for the membrane potential of a neuron has been extensively studied over the last decades. There seems to be a large agreement (cf.\ introduction Destexhe 1997) that the $4$-dimensional dynamical system proposed initially by Hodgkin and Huxley 1951 models adequately the mechanism of spike generation in response to an external input, in many types of neurons. It describes also the behavior of ion channels with respect to the two ion currents which are predominant: import of Na$^+$ and export of K$^+$ ions through the membrane, via voltage gated ion channels of a specific structure. For a modern introduction to the Hodgkin-Huxley model see Izhikevich 2009, in particular pp.\ 37--42 and figures 2.8 on  p.\ 33 and 1.7 on p.\ 5. 

The deterministic Hodgkin-Huxley system exhibits a broad range of possible and qualitatively quite different behavior of its solution, depending on the specific input. Desired periodic behavior (regular spiking of the neuron) appears only in special situations. See e.g.\  Rinzel and Miller 1980 for some interval $I$ such that time-constant input $c\in I$ results in periodic behavior of the solution; see e.g.\  Aihara, Matsumoto and Ikegaya 1984 for some interval $J$ such that oscillating input $t\to S(f\,t)$ with frequencies $f\in J$ (for some $1$-periodic function $S$) yields periodic behavior of the solution. In case of oscillating input, frequency has to be compatible with a range of 'preferred frequencies' of the Hodgkin-Huxley model, a fact which is similarly encountered in biological observations (see Izhikevich 2009, figure 1.7 on p.\ 5). There are also intervals $\wt I$ and $\wt J$ such that time-constant input $c\in \wt I$ or oscillating input at frequency $f\in \wt J$ leads to chaotic behavior of the solution. Periodic behavior includes that the period of the output can be a multiple of the period of the input. Using numerical methods, Endler 2012, Section 2, gives a complete tableau.

The first important question that one has to face when considering stochastic Hodgkin-Huxley systems is how to model the synaptic input received by the neuron from the other neurons. Actually, this question is a particular case of a more general problem, which is: in which way should 'noise' be included in a deterministic system such as Hodgkin-Huxley and what happens if one adds 'noise' to the system? 
There are some simulation studies (e.g., periodic signals embedded in Ornstein-Uhlenbeck type processes: Pankratova, Polovinkin and Mozekilde 2005, Yu, Wang, Wang and Liu 2001), but not much seems to be known mathematically. 

In our case, the input is driven by some deterministic $T-$periodic signal $S (t) $ which is randomly perturbed. We think of a cortical neuron which receives this input from its dendritic system. This dendritic system has a complicated topological structure and carries a large number of synapses which register spike trains emitted from a large number of other neurons within the same active network. There are statistical reasons to believe that the cumulated input as a function of time is well modeled by a time inhomogenous diffusion process $(\xi_t)_{t \geq 0}$ which is either of Ornstein-Uhlenbeck or of Cox-Ingersoll-Ross type, see H\"opfner 2007. More precisely, $(\xi_t)_{t \geq 0}$ is the strong solution to the SDE of mean-reverting type 
$$
d\xi_t \;=\; \left( S(t)-\xi_t \right) dt \;+\; \wt\si(\xi_t)\, dW_t 
$$
whose coefficients are such that 'periodic ergodicity' (cf.\ H\"opfner and Kutoyants 2010, section 2) holds for $\xi$. The signal $S$ is present e.g.\ in mean values $t\to E(\xi_t)$ through some deterministic transformation of $S$. The stochastic Hodgkin-Huxley model which we consider is then made of the four classical Hodgkin-Huxley variables together with $\xi,$ see equation ($\xi$HH) in Section \ref{sec:2.5} and also equation (HH) in Section \ref{sec:1.1}. It is therefore a $5$-dimensional SDE having the one-dimensional standard Brownian motion $W$ driving $\xi$ as the only source of 'noise'. For this reason it is a highly degenerate model in the sense that neither ellipticity nor the strong H\"ormander condition are fulfilled. Actually, only the weak H\"ormander conditions holds, and only locally.

Our model includes the feature of periodic behavior in the sense of a periodic structure of the semigroup of a Markov process. Several questions arise in this context: 
Does the noise have influence on the spiking behavior of the system? Does the noise enable a stochastic system to do what a comparable deterministic system would be unable to do? Moreover, from a probabilistic point of view it is also natural to determine whether there exist continuous transition densities for the system. This would show that the interaction between noise and drift can be strong enough to smoothen the degenerate $5$-dimensional diffusion.

Concerning the last issue, it is now classical to make use of Malliavin calculus techniques relying on the so-called H\"ormander condition. 
The H\"ormander condition is satisfied if some Lie algebra generated by coefficients of the system has sufficiently high dimension. For the strong H\"ormander condition, one makes only use of the diffusion coefficients to compute brackets, where-else in the weak case one can  also include the drift. 

In our case we only have the weak H\"ormander condition, and only locally, and our system is time inhomogeneous. In the non time homogeneous case, to the best of our knowledge, the existing results all require at least the strong H\"ormander condition, see Cattiaux and Mesnager 2002 and the references therein. Recently, Bally 2007, Fournier 2008 and De Marco 2012 considered the case of local ellipticity (in particular, the strong H\"ormander holds locally) in a time homogeneous framework with locally smooth coefficients. 
Using a technique based on estimates of the Fourier transform, introduced in these papers, we show that continuous transition densities indeed exist locally in all neighborhoods of points where the weak H\"ormander condition is satisfied. This is the content of Theorem \ref{theo:main} and Theorem \ref{theo:main2}. More precisely, we use a localization argument which is based on ideas of De Marco 2012;  technically our frame is more difficult since our system is not homogeneous in time and since we only have the weak H\"ormander condition locally. 

A natural question in this context is to exhibit an explicit set of points where the weak H\"ormander condition is satisfied. We can show numerically that on a specific segment of stability points for the deterministic Hodgkin Huxley system with constant input, the local H\"ormander condition is satisfied. Hence, locally at such stability points, continuous transition densities exist. We also can consider numerically a stable orbit of the deterministic system with constant input (sufficiently high), where a specific part of the  orbit --when the membrane potential up-crosses the resting level-- belongs to the set of points where the local H\"ormander condition is satisfied. However, in both cases, the weak H\"ormander condition 
neither is satisfied at all stability points, nor --and by far not-- at all points on the stable orbit. 

The only existence of continuous transition densities does not imply their strict positivity. Using a control argument, we prove in Theorem \ref{theo:pos} their strict positivity at stability points and at points on the stable orbit where the weak H\"ormander condition is satisfied.
Theorem \ref{theo:pos} is interesting also for the following reason:
%
it shows that with positive probability, our stochastic Hodgkin-Huxley system with $T$-periodic signal $S$ can imitate any deterministic Hodgkin-Huxley system driven by any $T$-periodic signal $\wt S$ over a certain time interval. Under some restriction on $\wt S$, this time interval can be arbitrarily long. More precisely, given a solution to the deterministic Hodgkin Huxley system associated to $\wt S$ over some time interval, small uniform tubes  around this deterministic solution will be in the support of the law of the stochastic system with $T$-periodic signal $S$. Hence, the stochastic system will be able to reproduce regularly spiking behavior during some period, followed by completely irregular behavior during some other period, followed by sticking to some equilibrium point during again some other period of time. 

This gives an answer to one of our questions: The stochastic system with signal $S$ can --with positive probability-- mimick deterministic systems with arbitrary  $\wt S$ over some time. Another question however is not answered by this assertion: what will be typical features of the path of the stochastic Hodgkin-Huxley system with $T$-periodic signal $S$ in the lon run? We know that the semigroup has a $T$-periodic structure, but neither this nor the preceding assertion allows to deduce what the system will do 'typically' when time tends to $\infty$. This is the question of determining whether the $5$-dimensional stochastic system is periodically ergodic which is outside the scope of the present paper.

Our paper is organized as follows. We present the deterministic and the stochastic Hodgkin-Huxley system in Section \ref{sec:main}. This section contains the main results, Theorem \ref{theo:hoer}, Theorem \ref{theo:main} and Theorem \ref{theo:pos}, on the existence of continuous transition densities and their positivity. Section \ref{sec:mall} is devoted to the study of smoothness properties of densities for strongly degenerate inhomogeneous SDE's and contains Theorem \ref{theo:main2} which is stated in a general frame, independently of the Hodgkin-Huxley model. The control argument is given in Section \ref{sec:th3}. The explicit calculation of the Lie brackets is postponed to Section \ref{section:th1}.

\section{Deterministic and stochastic Hodgkin-Huxley system. Main results.}\label{sec:main}

We will consider a neuron modeled by a Hodgkin-Huxley system which receives a periodic input $S$ from its dendritic system. The input is random and there are statistical reasons to believe that, as a function of time, this random input is well modeled by a time inhomogeneous diffusion of mean reverting type, see H\"opfner 2007.

We start by recalling briefly the deterministic model.

\subsection{HH with deterministic $T$-periodic input}\label{sec:1.1} 

Let a $T$-periodic deterministic signal $t\to S(t)$ be given. The Hodgkin-Huxley equations with input $S(t) dt$ are 
$$
\left\{\begin{array}{l}
dV_t  \;=\; S(t)\, dt \;- \left[\, \ov g_{\rm K}\,n_t^4\, (V_t-E_{\rm K}) \;+\; \ov g_{\rm Na}\,m_t^3\, h_t\, (V_t  -E_{\rm Na}) \;+\; \ov g_{\rm L}\, (V_t-E_{\rm L}) \right] dt\\
dn_t \;=\;  \left[\, \al_n(V_t)\,(1-n_t)  \;-\; \beta_n(V_t)\, n_t  \,\right] dt  \\
dm_t \;=\;  \left[\, \al_m(V_t)\,(1-m_t)  \;-\; \beta_m(V_t)\, m_t  \,\right] dt  \\
dh_t \;=\;  \left[\, \al_h(V_t)\,(1-h_t)  \;-\; \beta_h(V_t)\, h_t  \,\right] dt ,
\end{array}\right. 
\leqno{\rm (HH)}$$
where
$$\ov g_{\rm K} = 36, \; \ov g_{\rm Na} = 120, \; \ov g_{\rm L} = 0. 3 , \;E_{\rm K} = - 12, \; E_{\rm Na} = 120, \; E_{\rm L} =  10.6 ,$$ 
with notations and constants of Izhikevich 2009, pp.\ 37--38. The functions $\al_n, \beta_n, \al_m, \beta_m, \al_h, \beta_h$ in (HH) take values in $(0,\infty)$  and admit a power series representation on $\bbr$. They are given as follows.
\begin{equation}
\begin{array}{llllll}
\alpha_n(v)  &=& \frac{0.1-0.01v }{\exp(1-0.1v)-1}, &  \beta_n(v) &= &0.125\exp(-v/80) ,  \\
\alpha_m(v)& = &\frac{2.5-0.1v}{\exp(2.5-0.1v)-1} , & \beta_m(v)&= &4\exp(-v/18) ,  \\
\alpha_h(v) &= &0.07\exp(v/20) , &\beta_h(v) &=& \frac{1}{\exp(3-0.1v)+1}.
\end{array}
\end{equation}
Define for $v\in\bbr$
\begin{equation}\label{eq:ninfty}
n_\infty(v) := \frac{\al_n}{\al_n+\beta_n}(v) \;,\; m_\infty(v) := \frac{\al_m}{\al_m+\beta_m}(v) \;,\; h_\infty(v) := \frac{\al_h}{\al_h+\beta_h}(v) \;.
\end{equation}
If we think of keeping the variable $V$ constant in (HH), then these are equilibrium values in $(0,1)$ for the variables $n$, $m$, $h$ when $V\equiv v\in\bbr$.  

Write $\mathbb{Y}:= (V,n,m,h)$ for the 4d system of 'biological variables'. Sometimes we associate a fifth variable $J$ with $dJ_t = S(t) dt$ to the system and write $\mathbb{X} := (V,n,m,h, J)$ for the 5d system. Fixing some open interval $U \subset \bbr$ large enough to contain all values of the integrated signal $J,$ we write 
\beao
E_4 := \bbr\times (0,1)^3 \quad\mbox{state space of $\mathbb{Y}$, with points $(v,n,m,h)$} \,, \\
E_5 := \bbr\times (0,1)^3\times U \quad\mbox{state space of $\mathbb{X}$, with points $(v,n,m,h,\zeta)$}  
\eeao  
(see Section \ref{sec:4.4} for a proof of the fact that the system stays in $E_4$ whenever it starts there),
and use notation $F:E_4\to \bbr$ for drift terms not related to the signal in the first equation of (HH): 
\beam\label{eq:F}
F(v,n,m,h) &:=& \ov g_{\rm K}\,n^4\, (v-E_{\rm K}) \;+\; \ov g_{\rm Na}\,m^3\, h\, (v  -E_{\rm Na}) \;+\; \ov g_{\rm L}\, (v-E_{\rm L}) \nonumber \\
&:=& 36\,n^4\, (v+12) \;+120\,m^3\, h\, (v -120) \;+ 0.3\, (v-10.6) \;. 
\eeam

Define from (\ref{eq:F}) a function $F_\infty : \bbr\to\bbr$ by
\begin{equation}\label{eq:Finfty}
F_\infty(v) \;:=\;  F\left( v, n_\infty(v), m_\infty(v), h_\infty(v) \right) \;,\; v\in\bbr \;.
\end{equation}
In particular, if we select $c\in\bbr$ such that  $c = F_\infty(v),$ then
\begin{equation}\label{eq:equi}
( v, n_\infty(v), m_\infty(v), h_\infty(v) )\;\in\; E_4
\end{equation}
is an equilibrium point for the deterministic system (HH) with constant signal $S(\cdot)\equiv c  .$

\begin{ex}\label{ex:2}
It is well known that for sufficiently large values of constant signal $S(\cdot)=c$, the deterministic system (HH) exhibits regular spiking (see Rintzel and Miller 1980, for the model constants used here see Endler 2012, section 2.1, in particular figure 2.6). This means that for such values of $c$, the equilibrium point (\ref{eq:equi}) is unstable, and that there is a stable orbit for the $4$-dimensional system $\mathbb{Y}$ of 'biological variables'.
\end{ex}

\subsection{$T$-periodic diffusions carrying the signal $S$ and HH system with stochastic input $t\to \xi_t$ }\label{sec:2.5}

Take the $T$-periodic signal $t\to S(t)$ of subsection \ref{sec:1.1} and suppose moreover that $ t \to S(t) $ is smooth. Consider a diffusion 
\begin{equation}\label{eq:xi}
d\xi_t \;=\; (\, S(t)-\xi_t\,)\, \tau dt \;+\; \gamma\, q(\xi_t)\, \sqrt{\tau} dW_t 
\end{equation}
for suitable $q (\cdot)$, where we have chosen a parametrization in terms of $\tau$ (governing 'speed' of the diffusion) and $\gamma$ (governing 'spread' of one-dimensional marginals). We assume that the process $\left(\xi_t\right)_{t\ge 0}$ takes values in an open interval $U$ in $\bbr$, that $q (\cdot)$ is strictly positive on $U$, and that in restriction to every compact interval in $U,$ the function $q (\cdot)$ is of class $C^\infty ,$ bounded together with all derivatives. 
Then $\xi$ is a non time-homogenous diffusion which carries the signal $S$. We assume that $q (\cdot)$, $\tau$ and $\gamma$ are such that strong solutions to (\ref{eq:xi}) exist and such that the following holds: \\

{\bf (V1): } The grid chain $\left( \xi_{kT} \right)_k$ is positive Harris with invariant law $\mu$ on $U.$ \\

The $T$-periodic structure of the semigroup of transition probabilities of $(\xi_t)_{t\ge 0}$ combined with {\bf (V1)} implies, for arbitrary choice of a shift $0\le s<T$, that segment chains  
$$
\left(\, \xi_{[\,s+kT \,,\, s+(k+1)T \,]}\,\right)_{k\in\bbn_0} \quad\mbox{are positive Harris with invariant law on $C([0,T])$ denoted by  $m^{(s)}$} \;.
$$
It also implies that, for every $\ell\in\bbn$, $\ell$-segment chains 
$$
\left(\, \xi_{[\,s+k(\ell T) \,,\, s+(k+1)(\ell T) \,]}\,\right)_{k\in\bbn_0} \quad\mbox{are positive Harris with invariant law on $C([0,\ell T])$ denoted by  $m^{(s,\ell)}$} \;. 
$$
As a consequence, under {\bf (V1)}, trajectories of the process $\xi$ should in some sense get 'close' to the deterministic $T$-periodic signal $t\to S(t)$ as $t\to\infty$, for arbitrary choice of a starting point in $U$. In the next example we introduce two basic models that we have in mind: Cox-Intersoll-Ross and Ornstein-Uhlenbeck type $T-$periodic diffusions carrying the signal $S.$

\begin{ex}
a) CIR type: for some constant $K$ such that $K > \frac{\gamma^2}{2}+\sup|S|$, we take $U=(-K,\infty)$ and $q (x)=\sqrt{(x+K)\vee 0\;}$ for $x\in U$. By choice of the constant $K,$ the process $\xi$ will never attain $-K .$

We have Laplace transforms for $ \tilde \xi_t = \xi_t + K ,$ given $\tilde \xi_s = \tilde x > 0,$ which have the form 
$$ \lambda \to \exp \left\{ - \tilde x \Psi_{s,t} ( \lambda ) - \int_s^t \tilde S(v) \Psi_{v,t} ( \lambda ) \tau dv \right \},$$
where
$$ \Psi_{s,t} ( \lambda ) = \frac{\tau e^{ - \tau ( t-s) }}{ 1 + \lambda \frac{\gamma^2 }{2} ( 1 - e^{ - \tau ( t-s) } ) } , \; s < t , \; \lambda \in [0, \infty ) $$
and $ \tilde S (v) = S(v) + K .$  Note that 
$$ \Psi_{t_1, t_2 } \circ \Psi_{t_2, t_3} = \Psi_{ t_1, t_3 } \mbox{ on $[0, \infty ) ,$ for $ 0 \le t_1 < t_2 < t_3 < \infty ,$ } $$
and compare to the (time homogeneous) formulas (1.7)+(1.8), (1.12)+(1.13), (1.14) of Kawazu and Watanabe 1971.

If we write $\tilde \mu$ for the invariant law of $ (\tilde \xi_{kT} )_{k \in \bbn},$ $T-$periodicity of $\tilde S$ allows to write the Laplace transform of $\tilde \mu $ as 
$$
\la \;\;\lra\;\; \exp\left\{ -\int_{-\infty}^0 \tilde S(v)\, \Psi_{v,0}(\la)\, \tau dv \right\}  .
$$ 
Similarly, the invariant law of $(\tilde \xi_{kT +s } )_{ k \in \bbn} $ for $ 0 < s < T $ has Laplace transform 
$$ \la \;\;\lra\;\; \exp\left\{ -\int_{-\infty}^s \tilde S(v)\, \Psi_{v,0}(\la)\, \tau dv \right\}  .
$$ 
Taking derivatives in the last expression and noticing that $ (\frac{\partial}{\partial \lambda} \Psi_{v,t} )(0+) = e^{ - \tau (t-v)},$ 
expectations of $\tilde \xi_s $ starting at time $t=0 $ from $\tilde \xi_0 \sim \tilde \mu $ 
$$M(s):=  E_{ \tilde \mu , 0 } ( \tilde \xi_s)  = \int_0^\infty \tilde S ( s - \frac{r}{\tau} ) e^{ - r } dr $$
are $T-$periodic functions in $s.$

b) OU type: we take  $U=\bbr$ and $q (\cdot)\equiv 1$. Then we have an explicit representation 
$$ \xi_t = x e^{ - \tau ( t- s)} + \int_s^t e^{ - \tau (t-v) }
\left( \tau S(v) dv + \gamma \sqrt{\tau} d W_v \right) , \; t \geq s ,
$$
for the process starting at time $s$ in $x.$ With the same function $ s 
\to M(s) = \int_0^\infty S( s - \frac{r}{\tau} ) e^{ - r } dr $ as in a), the invariant law $\mu $ of $(\xi_{kT})_{k \in \bbn} $ is 
$$ \mu = {\cal N} ( M(0), \frac{\gamma^2 }{2} ) $$
and the law of $\xi_s $ starting at time $t=0$ from $\xi_0 \sim \mu$ is 
$$ {\cal L}_{ \mu , 0} ( \xi_s) = {\cal N} ( M(s) , \frac{ \gamma^2 }{2} ) .$$
(cf. H\"opfner and Kutoyants 2010, Ex. 2.3). 

Hence in both cases a) and b), the $T-$periodic signal $S(\cdot) $ is 
expressed in the process $\xi $ under 'periodically invariant' regime in form 
of moving averages 
$$ s \to E_{\mu, 0} ( \xi_s ) = M(s) = \int_0^\infty S(s - \frac{r}{\tau} ) e^{ - r } dr $$
which are $T-$periodic.  
\end{ex}

Consider now the HH equations driven by stochastic input $d\xi_t$, i.e.\ the 5d system 
$$
\left\{\begin{array}{l}
dV_t  \;=\; d\xi_t \;- \left[\, \ov g_{\rm K}\,n_t^4\, (V_t-E_{\rm K}) \;+\; \ov g_{\rm Na}\,m_t^3\, h_t\, (V_t  -E_{\rm Na}) \;+\; \ov g_{\rm L}\, (V_t-E_{\rm L}) \right] dt\\
dn_t \;=\;  \left[\, \al_n(V_t)\,(1-n_t)  \;-\; \beta_n(V_t)\, n_t  \,\right] dt  \\
dm_t \;=\;  \left[\, \al_m(V_t)\,(1-m_t)  \;-\; \beta_m(V_t)\, m_t  \,\right] dt  \\
dh_t \;=\;  \left[\, \al_h(V_t)\,(1-h_t)  \;-\; \beta_h(V_t)\, h_t  \,\right] dt  \\
d\xi_t \;=\; (\, S(t)-\xi_t\,)\, \tau dt \;+\; \gamma\, q (\xi_t)\, \sqrt{\tau} dW_t 
\end{array}\right.
\leqno{\rm (\xi HH)}
$$
under assumption {\bf (V1)}. Write $E_5=\bbr\times (0,1)^3\times U$ for the corresponding state space and denote its elements by $(v,n,m,h,\zeta)$. Let $\left(\,P_{s_1,s_2}(x_1, dx_2)\,\right)_{0\le s_1<s_2<\infty}$ denote the semigroup of transition probabilities (which is non-homogenous in time) of the 5d system ${\rm (\xi HH)}$. Due to $T$-periodicity of the deterministic signal $t\to S(t)$, the semigroup is $T$-periodic in the following sense:  
$$
P_{s_1,s_2}(x_1, dx_2) \;=\; P_{s_1+kT,s_2+kT}(x_1, dx_2) \quad\mbox{for all $k\in\bbn_0$} \;. 
$$

\subsection{Existence of densities for the stochastic HH system}
In order to state our first theorem we have to introduce some notation. Let us first denote the drift terms related to $n,m,h$ in (HH)  
\beao 
G_n(v,n) =  \al_n(v)\,(1-n)  \;-\; \beta_n(v)\, n &,&  g_n(v,n) =  \al'_n(v)\,(1-n)  \;-\; \beta'_n(v)\, n \,,\\
G_m(v,m) =  \al_m(v)\,(1-m)  \;-\; \beta_m(v)\, m &,&  g_m(v,m) =  \al'_m(v)\,(1-m)  \;-\; \beta'_m(v)\, m \,, \\
G_n(v,h) =  \al_h(v)\,(1-h)  \;-\; \beta_h(v)\, h &,&  g_h(v,h) =  \al'_h(v)\,(1-h)  \;-\; \beta'_h(v)\, h \,  ,
\eeao
where notation $'\,$ is reserved for derivative with respect to $v.$ Then let   
\begin{equation}\label{eq:det}
D(v,n,m,h) \;:=\; \det \left( \begin{array}{lll} 
g_n' &  g_n'' &  g_n'''  \\
g_m' &  g_m'' &  g_m'''  \\
g_h' &  g_h'' &  g_h'''  \\
\end{array} \right) (v,n,m,h) \quad,\quad (v,n,m,h)\in E_4 \;. 
\end{equation}
This determinant will be important in the sequel. We introduce
$$ {\cal O} := \{ (v,n,m,h) \in E_4 : D( v,n,m,h) \neq 0 \} .$$

Notice that by continuity of $D$ on $E_4,$ the set $ {\cal O} $ is open. Moreover, $\lambda ( {\cal O}^c) = 0 .$ This can be seen as follows. Firstly it can be shown numerically that ${\cal O} $ is not empty (see Section \ref{ex:3} below). Moreover, for any fixed $v \in \bbr ,$ the function $ (n,m,h) \to D( v, n,m,h) $ is a polynomial of degree three in the three variables $n,m,h .$ In particular, for any fixed $v,$ either $D(v,.,.,.)$ vanishes identically on $(0,1)^3$, or the zeros of $(n,m,h) \to  D(v,n,m,h)$ form a two-dimensional sub-manifold of $(0,1)^3 .$ 
Finally, since the determinant is a sum of terms
$$
\mbox{(some power series in $v$)} \cdot n^{\vep_n} m^{\vep_m} h^{\vep_h}
$$
with epsilons taking values 0 or 1, it is impossible to have small open $v$-intervals where $D(v,.,.,.)$ vanishes identically on $(0,1)^3 .$ Then integrating over $v$, Fubini finishes the proof.

Now we have the following result.

\begin{theo}\label{theo:hoer}
The weak H\"ormander condition holds at points $x=(v,n,m,h,\zeta)$ in $E_5$ whose  first four components belong to ${ \cal O} . $ 
\end{theo}

We provide in Section \ref{ex:3} below a numerical study of the set ${\cal O} $ where the H\"ormander condition holds. 
Once the H\"ormander condition holds locally, we are able to show that the process, in spite of its very degenerate structure (only the first and the fifth variable carry Brownian noise), possesses Lebesgue densities locally. This is the content of the next theorem.

\begin{theo}\label{theo:main}
For $0\le s_1<s_2<\infty$, consider the 5d process ${\rm (\xi HH)}$ starting at time $s_1\ge 0$ from arbitrary $x\in E_5$. Then in restriction to the subset ${\cal O} {\times}U$ 
of $E_5\,$, the law $\,P_{s_1,s_2}(x, \cdot)\,$ admits a continuous Lebesgue density $\,p_{s_1,s_2}(x, \cdot)\,$. Moreover, for any fixed $ x' \in {\cal O} \times U , $ the map $ x \to p_{s_1,s_2}(x, x') $ is lower semi-continuous. 
\end{theo}

Note that this is a local result in the second variable for fixed starting point, local in restriction to ${\cal O} \times U .$ 
In particular we impose the H\"ormander condition on the second variable and not on the starting point.

\subsection{Numerical study of the determinant $D$}\label{ex:3}
We study numerically the above determinant (\ref{eq:det}) and provide some figures. First, we can not expect to have $D(v,\cdot,\cdot,\cdot)\ne 0$ on $(0,1)^3$. Indeed, we find that $D$ vanishes at points $ (v, n_\infty (v), m_\infty (v) , h_\infty (v))$ with 
$v$ located at $\approx -11.4796$ and $\approx +10.3444$.
Second, $\cal O$ is certainly non-empty since we find a strictly negative value of the determinant e.g. at the point $(0, n_\infty(0), m_\infty(0), h_\infty(0) )$.

In order to obtain more detailed information about ${\cal O},$ we calculate the determinant $D$ in stable equilibrium points and along stable orbits of (HH) corresponding to different constant inputs $S(\cdot) \equiv c. $ 

First of all, we find that the function 
$
v \;\lra\; D\left(v, n_\infty(v), m_\infty(v), h_\infty(v)\right)
$
has exactly two zeros on $I_0 := (- 15, + 30)$ which are located at $v \approx -11.4796$ and $\approx +10.3444$.
Since the function $F_\infty (v) $ of (\ref{eq:Finfty}) is strictly increasing on a large interval containing $I_0,$ all points $(v, n_\infty(v), m_\infty(v), h_\infty(v))$ with $v \in I_0$ correspond to equilibrium points of (HH) associated to constant input $c$ where $ F_\infty (v) = c .$ The corresponding range of values for $c$ is given by $c \in  (F_\infty ( -10) , F_\infty  (+10)) = (-6.15, 26.61) .$ For $v \approx 0$ we find $c\approx 0.0534 .$
Hence for all values of $c$ belonging to $  (-6.15, 26.61),$ the determinant of the associated equilibrium point stays strictly negative. 

Moreover, also on stable orbits of (HH) with large enough constant signal, we can not expect that $D $ never vanishes. Indeed (see Figure 1 below), on a good approximation to the stable orbit for constant signal $S(\cdot)=15$, we find numerically a segment requiring approximately one third of the time needed to run the orbit (very roughly, this segment starts when the variable $v$ up-crosses the level $-2$ and ends when $v$ up-crosses the level $+5$) where the determinant  $D(v,n,m,h)$ in (\ref{eq:det}) is negative and well separated from zero. On the remaining parts of the orbit, the determinant changes sign several times, and in particular takes values very close to zero immediately after 'the top of the spike', i.e.\ after the variable $v$ has attained its maximum over the stable orbit.

In Figure 2 below we consider a deterministic HH with constant input $c=15$ starting in a numerical approximation to its equilibrium point. It is seen that this equilibrium point is unstable, and the system switches towards a stable orbit. In this picture, already the last four orbits can be superposed almost perfectly. Figure 1 shows the value of the determinant at equidistant time epochs on the last complete orbit (starting and ending when the membrane potential $v$ up-crosses the level $0$, and having its spike near time $t=180$). 

\newpage

\begin{figure}[!h]
\begin{center}
   \includegraphics[width=0.80\textwidth]{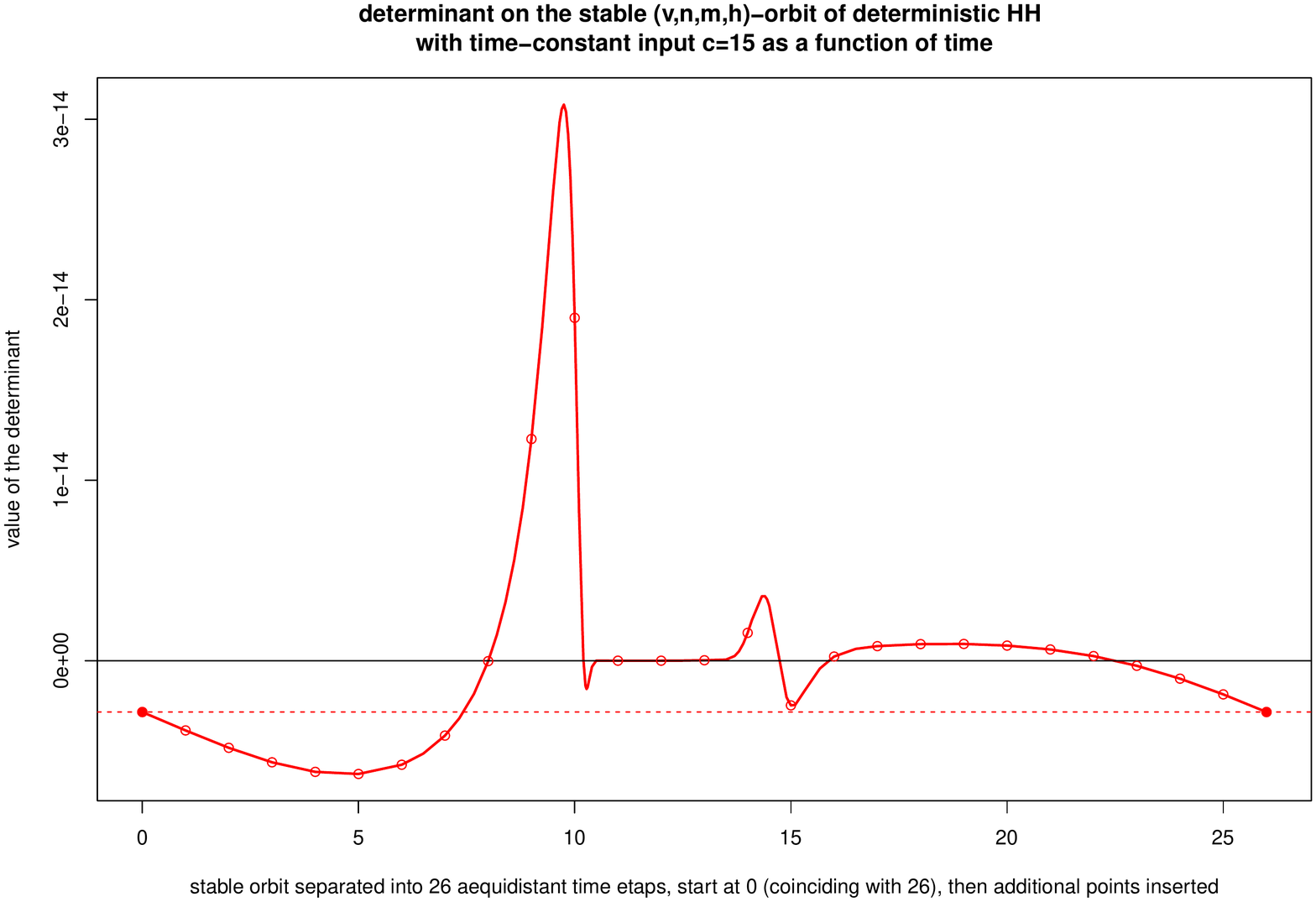} 
\end{center}
\caption{\small Determinant $D$ of (\ref{eq:det}) calculated on the stable orbit $t \to (v_t,n_t,m_t,h_t)$ of the deterministic system (HH) with constant input $c=15$. The time needed to run the stable orbit is $\approx 12.56$ ms.}
\begin{center}
   \includegraphics[width=0.8\textwidth]{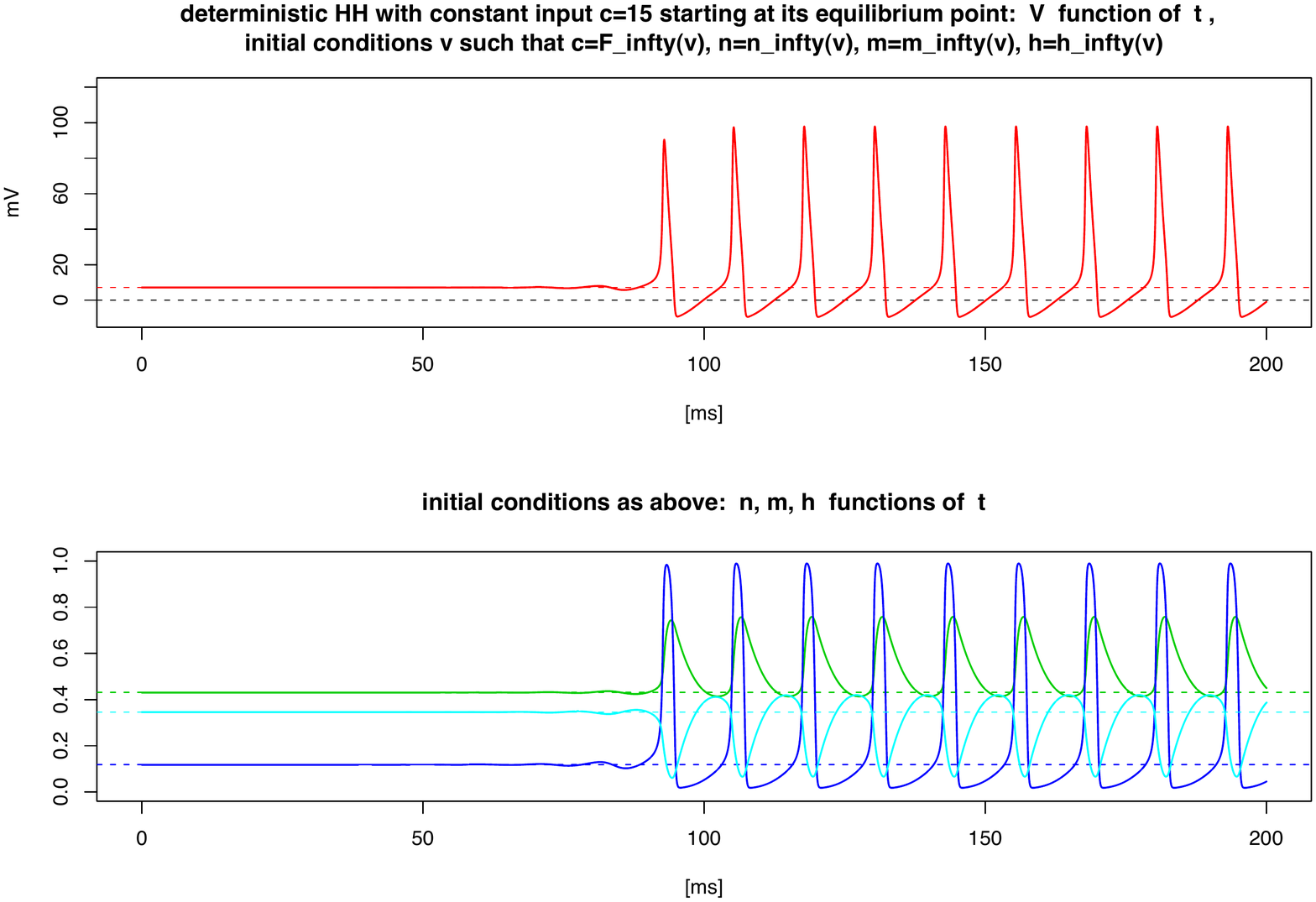} 
\end{center}
\caption{\small Deterministic HH with constant input $c=15.$
}
\end{figure}

\subsection{The HH system with stochastic input can reproduce any deterministic feature}
In our Theorem \ref{theo:main}, no condition is imposed on the starting point, and the density might be identically $0$ on ${\cal O} \times U$. In order to exhibit regions of the space where the transition density is strictly positive we use control arguments. These arguments are interesting also for the following reason. They show that the stochastic system ($\xi $HH) is able to reproduce any possible deterministic spiking behavior during any arbitrary long period, in the sense that any of these behaviors is in the support of the law of the process. In other words, the stochastic system ($\xi$HH) driven by the signal $S$ can stay with positive probability for an arbitrarily long time in arbitrarily small tubes around deterministic solutions of systems (HH) driven by any smooth $T$-periodic signal.

During this subsection, we still denote $ t \to S(t)$ the fixed $T-$periodic signal which is carried by the diffusion process $\xi_t $ and which governs the evolution of ${\rm (\xi HH)} .$ Moreover, $t  \to \tilde S (t)$ will denote any other signal chosen independently of $S .$ We shall write $\mathbb{P}_{0, x} $ for the law of the process ${\rm (\xi HH)} $ starting from $ x \in E_5 $ at time $0.$ 

\begin{theo}\label{theo:pos}
Fix $ (v,n,m,h ) \in E_4 $ and $t > 0 .$ Then for any smooth $T-$periodic signal $s \to \tilde S(s)$ and any initial value 
$ \zeta \in U$ such that $ \tilde J_s := \zeta + \int_0^s \tilde S (u) du \in U $ for all $ s \le t ,$ the following holds:
 
Let $\mathbb{X}_s, s \le t ,$ be the associated deterministic system (HH), driven by the signal $ s \to \tilde S (s).$ 
Write $ x = ( v,n,m,h, \zeta ) .$ Then we have for any $\varepsilon > 0 , $ 
$$ \mathbb{P}_{ 0,x} \left( \left\{ f \in C ( [0, \infty [ , \bbr^5 ) : \sup_{ s \le t } | f(s) - \mathbb{X}_s | \le \varepsilon \right\} \right) > 0 .$$
Moreover, for any fixed $ \varepsilon > 0,$ there exists $\delta > 0 ,$ such that for all $ x'' \in B_\delta ( x) , $ 
$$ \mathbb{P}_{ 0,x''} \left( \left\{ f \in C ( [0, \infty [ , \bbr^5 ) : \sup_{ s \le t } | f(s) - \mathbb{X}_s | \le \varepsilon \right\} \right) > 0 .$$
\end{theo}

We sketch two situations where the above theorem can be applied successfully. Consider first the situation of constant signal $ \tilde S \equiv c .$  
In what follows, $ B_{\varepsilon }( x)$ denotes the open ball of radius $\varepsilon$ centered in $x .$ Moreover, for a suitable choice of $c,$ let $v_c$ be such $F_\infty (v_c) = c.$  

\begin{cor}\label{cor:1}
Fix a constant $c  $ such that $v_c$ exists and let $ \zeta \in U$ such that $\tilde  J_s = \zeta + c s \in U $ for all $ s \le t  .$ Then for all $ \varepsilon > 0 ,$
$$ P_{ 0, t} ( x_c, B_{\varepsilon }( x_c') ) > 0  ,$$
where 
$$ x_c = ( v_c, n_\infty ( v_c), m_\infty (v_c), h_\infty ( v_c) , \zeta ) , \;  x_c' = (v_c, n_\infty ( v_c), m_\infty (v_c), h_\infty ( v_c)  , \zeta + ct ) .$$ 
Moreover, for any fixed $ \varepsilon > 0,$ there exists $\delta > 0 ,$ such that for all $ x'' \in B_\delta ( x_c) , $ 
$$    P_{ 0, t} ( x'' , B_{\varepsilon }( x_c') ) > 0  .$$
\end{cor} 

Now we combine this result with Theorem \ref{theo:main} above and use the fact that H\"ormanders condition holds for several stability points (cf. Section \ref{ex:3} above). 
We keep the notation of Corollary \ref{cor:1}.

\begin{cor}\label{cor:2}
For any $c$ such that $v_c $ exists and such that the local H\"ormander condition holds at $x_c,$ there exists
$ \delta > 0 $ such that for $ K_c = \overline B_\delta ( x_c)   $ and $ K_c' = \overline B_\delta (x_c') ,$  
$$ \inf_{ x \in K_c } \inf_{x' \in K_c'}  p_{ 0, t} (x, x') > 0 .$$
\end{cor}
The above result holds in particular for $-6.15 < c < 26.61$. 

The second situation which we consider is the deterministic system (HH) with sinusoidal signal 
$
\tilde S(t)  = a \left( 1 + sin( 2\pi \frac{t}{T} ) \right) , 
$
$a>0$ some constant. This system presents additional features (see Aihara, Matsumoto and Ikegaya 1984 for a modified system, for the (HH) system as above see Endler 2012 Ch.\  2.2). There are specified subsets $D_1$, $D_2$, $D_3$, $D_4$ in $(0,\infty)\times(0,\infty)$ with the following properties: i) for $(a,T)$ in $D_1$, the system (HH) is periodic with small oscillations which can not be interpreted as 'spiking'; ii) for  $(a,T)$ in $D_2$, the system moves on a $T$-periodic orbit, and the projection $t\to V_t$ resembles the membrane potential of a regularly spiking neuron (single spikes or spike bursts per orbit); iii) for  $(a,T)$ in $D_3$, the system moves on a $kT$-periodic orbit for some multiple $k\in\bbn$; iv) for $(a,T)$ in $D_4$, the system behaves 'irregularly' and does not exhibit periodic behavior. 

\begin{cor}\label{cor:3}
Let $ (a , T) \in D_2 $ and $ ( 0, n^* , m^* , h^* ) $ be a point on the $T-$periodic orbit of the associated deterministic system (HH) such that   
\begin{equation}
 ( 0, n^* , m^* , h^* ) \in {\cal O} .
\end{equation}
Fix $ \zeta \in U $ such that $ \zeta + \int_0^s \tilde S(u) du \in U $ for all $ 0 \le s \le T$ and write $ x^* = ( 0, n^* , m^* , h^* , \zeta ) $ and $ z^* =( 0, n^* , m^* , h^* , \zeta + \int_0^T \tilde S(u) du ) .$ Then there exists
$ \delta > 0 $ such that for $ K  = \overline B_\delta ( x^*)   $ and $ K' = \overline B_\delta (z^*) ,$  
$$ \inf_{ x \in K  } \inf_{x' \in K' }  p_{ 0, T} (x, x') > 0 .$$
\end{cor}

The proofs of Theorem \ref{theo:hoer}, Theorem \ref{theo:main} and Theorem \ref{theo:pos} are given in Sections \ref{section:th1}, \ref{sec:mall} and \ref{sec:th3}, respectively.

\section{Discussing $T-$periodic ergodicity of the stochastic (HH) system} 
Using Lyapunov functions, one can show that the five dimensional diffusion process $X$ of ($\xi$HH) is ultimately bounded (cf. Miyahara 1972). As a consequence, the grid chain $ ( X_{kT})_{ k  } $ possesses (possibly infinitely many) invariant probability measures. 

Suppose there exists some orbit and an associated recurrent point $x^* $ as in Corollary \ref{cor:3} above such that, for all starting points, the grid chain visits $K  = \overline B_\delta ( x^*) $ infinitely often. Then Corollary \ref{cor:3} implies the Harris property of the grid chain, and hence $T-$periodic ergodicity of the stochastic Hodgkin-Huxley model (compare to H\"opfner and L\"ocherbach 2011), the sets $ \overline B_\delta ( x^*) $ playing the role of 'small' sets of the system. 

Hence Corollaries \ref{cor:2} and \ref{cor:3} above are important steps towards periodic ergodicity of the process. However, in order to show the recurrence property of one of these sets $K,$ one has to find a Lyapunov function which forces the system to come back precisely to this set $K.$ This requires a more detailed study of the recurrence properties of the process and of the possible sets $K$ and will be part of some future work.

\section{Smoothness of densities of a strongly degenerate inhomogeneous SDE with locally smooth coefficients: proof of Theorem \ref{theo:main}}\label{sec:mall}
Let $X_t = (V_t, n_t, m_t, h_t, \xi_t )$ be the $5-$dimensional diffusion process of ($\xi \! \!$ HH). The aim of this section is to show that this process admits locally a continuous Lebesgue density. Classically, the main tool to
prove that the law of a diffusion admits a smooth density is Malliavin calculus. A usual technical condition is to suppose that the coefficients of the SDE are bounded $C^\infty -$functions with bounded derivatives of any order. This condition is obviously not satisfied in our situation. Moreover, in our case, a one-dimensional Brownian motion is driving a 5-dimensional system; as we will see in Section \ref{section:th1}, the H\"ormander condition holds only locally. Finally, the drift coefficient depends on time. 

Therefore we have to apply local results which are similar to those obtained by Kusuoka and Stroock 1985. The results obtained there hold only in a time homogeneous frame. So in what follows we extend the results of Kusuoka and Stroock to the non time homogeneous case. In order to do so, we recall ideas of De Marco 2011 and adopt them to our frame. The results we obtain are interesting in their own right, therefore we state them in a general setting.  

We start by introducing some notation and the general framework in which we will work.

\subsection{Notation}
Let $ m \geq 1 .$ We consider processes taking values in $\bbr^m$ and write $ x = (x^1, \ldots , x^m) $ for generic elements of $\bbr^m .$ We will identify the time variable $t$ with $x^0 .$ Let $\sigma $ be a measurable function from $\bbr^m $ to $\bbr^m $ and $b$ a a smooth function from $ [0, \infty [ \times \bbr^m $ to $\bbr^m . $ For $y_0 \in \bbr^m $ and $\delta > 0 $ we denote by $B_\delta (y_0) $ the open ball of radius $\delta $ centered in $y_0 .$ For any  open subset $ A \subset \bbr^m ,$ $C_b^\infty (A)$ denotes the class of infinitely differentiable functions defined on $A$ which are bounded together with all partial derivatives of any order. Fix some $ 0 < R \le 1. $ We consider the SDE 
\begin{equation}\label{eq:equationms}
X_t^i = x^i + \int_0^t b^i ( s, X_s) ds + \int_0^t \sigma^i (X_s) d W_s , \; t \geq 0 , \; i = 1, \ldots , m , 
\end{equation}
for all $ x \in \bbr^m .$ Here, $W$ is a one-dimensional Brownian motion and $ \sigma : \bbr^m \to \bbr^m $ is identified with an $m \times 1$matrix. 

We impose the following conditions on the coefficients of the above equation. 
$$ \mbox{ Existence of strong solutions holds for the couple $ (b(t,x) , \sigma (x))$ .} \leqno{\rm \bf (H1)} $$
Let then $ ( X_t, t \geq 0 )$ be a strong solution of (\ref{eq:equationms}). We suppose moreover that 
there exists a growing sequence of compacts $ K_n = [a_n, b_n] \subset \bbr^m  , \; K_n \subset K_{n+1} ,$ such that the following holds. If the starting point $x$ satisfies
$ x \in \bigcup_{n } K_n ,$ then we have 
$$ T_n := \inf \{ t : X_t \notin K_n \}  \to \infty  \mbox{ almost surely as } n \to \infty .$$
In the above, $ [a_n, b_n] = \prod_{i = 1}^m [a_n^i , b_n^i ] ,$ where $ a_n = ( a_n^1 , \ldots , a_n^m ) .$ Due to the above condition, we can introduce 
\begin{equation}\label{eq:e}
 E = \bigcup_{n } K_n ,
\end{equation}
which is the state space of the process.

We impose local smoothness on each compact $K_n .$ For that sake, fix some $T > 0 $ and suppose: For all $n,$ for all multi-indices $\beta ,$
we have
$$ \sigma \in C_b^\infty (K_n , \bbr^m ), \quad  b(t,x) +  \partial_{\beta}  b (t,x)  \mbox{ is bounded on  }  [ 0, T ] \times K_n , \leqno{\rm \bf (H2)}\\
$$
where for $\beta \in \{ 0, \ldots, m \}^l, l \geq 1,$  $\partial_\beta =\frac{ \partial^l}{ \partial x^{\beta_1} \ldots \partial x^{\beta_l} } . $ Recall that we identify $x^0 $ with $t.$  

\subsection{H\"ormander condition}\label{sec:hoerm} 
Due to the strong degeneracy of our biological system ($\xi \! \!$ HH), the condition of
ellipticity is no where satisfied. However, as we will see in Section \ref{section:th1} below,  
the H\"ormander condition holds locally. In order to state the H\"ormander condition, we have to rewrite the above equation (\ref{eq:equationms})
in the Stratonovitch sense. That means, we replace the drift function $ b(t,x) $ by $\tilde b (t,x) $ defined as
$$ \tilde b^i (t,x) = b^i (t,x) - \frac12 \sum_{k=1}^m \sigma^k (x) \frac{ \partial \sigma^i }{\partial x^k } (x) ,\; x \in E , \; 1 \le i \le m .$$
The above drift function is non-homogeneous in time.  The associated directional derivative is 
$$ V_0 = \frac{\partial }{\partial t } + \sum_{ i = 1}^m \tilde b^i (t,x) \frac{\partial }{\partial x^i } = \frac{\partial}{\partial t} + \tilde b .$$
Notice that $V_0$ can be identified with the $ (m+1 )-$dimensional function $ V_0 ( t, x ) = ( 1, \tilde b^1 , \ldots , \tilde b^m) .$

By convention, all other functions $ V (t,x) : [0, \infty [ \times \bbr^m \to \bbr^m $ different from $V_0$ will be interpreted only as directional derivatives with respect to space variables
$$ V(t,x) = \sum_{i=1}^m V^i (t,x) \frac{\partial}{\partial x^i},$$
even if they are time dependent.
Hence we identify $V(t,x) $ with the $(m+1)-$dimensional function $V(t,x) = (0, V^1 , \ldots , V^m ).$ 

Now we can introduce the successive Lie brackets. We start by putting $ V_1 (x) = \sigma (x) $ and identify this function with the directional derivative $ \sum_{ i =1}^m V_1^i (x) \frac{\partial }{ \partial x^i } = \sum_{ i =1}^m \sigma^i (x) \frac{\partial }{ \partial x^i }.$

We adopt the formalism of Kusuoka and Stroock 1985 and put $ A = \emptyset \cup \bigcup_{ l = 1}^\infty ( \{ 0 , 1 \})^l .$ For any $ \alpha \in A, $ define $ | \alpha | := l $ if $ \alpha \in \{ 0, 1\}^l, l \geq 1 , $ $ | \emptyset | = 0 .$ Moreover, let $ \| \emptyset \| = 0 $ and $ \| \alpha \| = | \alpha | + card \{ j : \alpha_j = 0 \} , $ if $ | \alpha | \geq 1 .$ Finally, we put $ \alpha ' = ( \alpha_1, \ldots , \alpha_{l-1}) $ if $\alpha = ( \alpha_1 , \ldots , \alpha_l ), l \geq 2 , $ $\alpha' = \emptyset $ if $l = 1 .$ 

Recall that $ V_0 = \frac{ \partial }{\partial t} + \tilde b = \frac{ \partial}{\partial x^0} + \tilde b $ and that $ V_1 = \sigma .$ For any $ V  : [ 0, \infty [ \times \bbr^m \to \bbr^m ,$ define inductively in $| \alpha |,$ 
$$ V_\emptyset (t,x) := V(t,x) $$
and for $| \alpha | \geq 1,$ 
\begin{equation}\label{eq:valpha}
 V_{(\alpha ) } ( t,x) := 
[ V_{ \alpha_l } , V_{( \alpha ')} ] 
.
\end{equation}
Here, $[V, W]$ denotes the Lie bracket defined by
$$ [V, W ]^i = \sum_{j = 0}^m \left( V^j \frac{ \partial W^i }{\partial x^j } - W^j \frac{ \partial V^i }{\partial x^j } \right) .$$
In other words, if $ V = V_1 = \sigma ,$ we have 
$$ [V_1, W]^i = \sum_{j = 1}^m \left( \sigma^j \frac{ \partial W^i }{\partial x^j } - W^j \frac{ \partial \sigma^i }{\partial x^j } \right) ,$$
and the time variable does not play any role. But if $ V = V_0, $ we have, since $ V_0^0 \equiv 1 ,$ 
$$ [V_0, W ]^i = \sum_{j = 0}^m \left( V_0^j \frac{ \partial W^i }{\partial x^j } - W^j \frac{ \partial V_0^i }{\partial x^j } \right) = 
\frac{ \partial W^i }{\partial t } +\sum_{j = 1}^m \left( V_0^j \frac{ \partial W^i }{\partial x^j } - W^j \frac{ \partial V_0^i }{\partial x^j } \right) .$$

Finally, for any $x \in E$ and any $\eta \in \bbr^m,$  we define 
$$ {\cal V}_L (t,x, \eta) = \sum_{ \alpha : \| \alpha \| \le L- 1 } < (V_1)_{(\alpha)} (t,x), \eta >^2  $$
and 
\begin{equation}\label{eq:lie1}
 {\cal V}_L (x) = \inf_{ 0 \le t \le T , \eta : \| \eta \| = 1 } {\cal V}_L (t,x, \eta ) \wedge 1 .
\end{equation}

We assume:
\begin{ass}
There exists $y_0  $ with $B_{5R} (y_0) \subset E$ and some $L \geq 1 $ such that the following local H\"ormander condition holds:
$$  \mbox{ We have } {\cal V}_L (y) \geq c( y_0, R)  > 0 \; \mbox{ for all } y \in B_{ 3 R} (y_0 ) . \leqno{\rm \bf (H3)}  $$
\end{ass}

Now our result is as follows.

\begin{theo}\label{theo:main2}
Assume {\bf (H1), (H2) } and {\bf (H3)}. Then for any initial condition $x \in E 
$ and for any $t \le T ,$
the random variable $X_t $ admits a Lebesgue density $ p_{ 0, t} ( x, y ) $ on $ B_R ( y_0 ) $ which is continuous with respect to $ y \in B_R ( y_0 )  .$ Moreover, for any fixed $ y \in B_R ( y_0) ,$ $E \ni  x \to  p_{0, t } (x,y) $ is lower semi-continuous.
\end{theo}

The proof of Theorem \ref{theo:main2} is given in the next subsection and uses localization arguments.

\subsection{Proof of Theorem \ref{theo:main2}} 
Recall that a random vector taking values in $\bbr^m $ is said to admit a density on an open set $O \subset \bbr^m $ if
$$ E (f(X)) = \int f (x) p(x) dx ,$$
for any continuous function $f \in C_b ( \bbr^m ) $ such that $ supp (f) \subset O ,$ for some positive function $ p \in L^1 (O) .$  We rely on the following classical criterion for smoothness of laws which is based on a Fourier transform method. 

{\bf 4.3 Proposition:} Let $\nu$ be a probability law on $\bbr^m $ and let $ \hat \nu (\xi)$ be its Fourier transform. If $\hat \nu $ is integrable, then $\nu$ is absolutely continuous and 
$$ p(y) = \frac{1}{(2 \pi)^m} \int_{\bbr^m} e^{ - i < \xi , y>} \hat \nu (\xi ) d \xi $$
is a continuous version of its density.  

We have to replace the above argument by a localized one. This localization follows ideas that have been developed by De Marco 2011 and that we adopt to our frame. We start by taking a function $\Phi \in C_b^\infty ( \bbr^m ) $ such that $ 1_{ B_R (0)} \le \Phi \le 1_{ B_{2R}(0)} .$ Fix $x$ and $t \le T$ and suppose that $ E_x ( \Phi ( X_t - y_0 )) := m_0 > 0 .$ Then we can define a probability measure $\nu$ via 
\begin{equation}\label{eq:mu}
 \int f(y) \nu (dy) : = \frac{1}{m_0} E_x \left( f( X_t) \Phi ( X_t - y_0 )\right) .
\end{equation}
In order to prove Theorem \ref{theo:main2} it is sufficient to show that $\nu$ admits a continuous Lebesgue density. For that sake let 
$$ \hat \nu (\xi ) = \frac{1}{m_0} E_x \left( e^{i < \xi , X_t>} \Phi ( X_t - y_0 ) \right) $$
be its Fourier transform. All we have to show is that $\hat \nu (\xi ) $ is integrable. In order to do so, we use Malliavin calculus localized around $y_0.$ More precisely, we localize the coefficients of the SDE (\ref{eq:equationms}) in the following way. Let $\psi \in C^\infty_b (\bbr^m ) $ such that 
$$ \psi (y ) = \left\{ 
\begin{array}{ll} 
y & \mbox{ if } |y| \le 4 R \\
5 R \frac{ y}{|y|} & \mbox{ if } |y| \ge 5 R 
\end{array}
\right. $$
and $|\psi (y)| \le 5 R $ for all $ y .$ Put $ \bar b (y) = b ( \psi ( y - y_0 ))$ and $\bar \sigma (y) = \sigma ( \psi ( y - y_0 )) .$ Then by condition {\bf (H2)}, $\bar b$ and $\bar \sigma$ are $C_b^\infty-$extensions of $ b_{| B_{4R} (y_0)} $ and $\sigma_{ | B_{4R} (y_0)} .$ 

We denote $ \XXX $ the unique strong solution of the equation 
\begin{equation}\label{eq:processgood}
\bar X_s^i = x^i + \int_0^s \bar b^i ( u, \bar X_u) du + \int_0^s \bar \sigma^i (\bar X_u) d W_u ,\;  u \le T ,\; 1 \le i \le m .
\end{equation}
Up to the first exit time of $ B_{4R} (y_0) ,$ both processes $ \XXX $ and $ X$ coincide. 

Now for some fixed $ \delta \in ] 0, t/2 \wedge 1[,$ we put 
$$ \tau_1 = \inf \{ s \geq t - \delta : X_s \in B_{3R} (y_0 ) \} \quad 
\mbox{
and }
\quad 
\tau_2 = \inf \{ s \geq \tau_1 : X_s \notin B_{ 4 R }(y_0) \} .$$  
Then, 
\begin{multline*}
\{ \Phi ( X_t - y_0 ) > 0 \} = \{ \Phi ( X_t - y_0 ) > 0 ; t - \delta = \tau_1 < t < \tau_2  \} \\
\cup \left\{ \Phi ( X_t - y_0 ) > 0 ; \sup_{ 0 \le s \le \delta } | \XXX_s ( X_{\tau_1}) - X_{ \tau_1} | \geq R \right\} .
\end{multline*}
Hence,
\begin{multline}\label{eq:fourier}
m_0 \hat \nu ( \xi ) = E_x \left( e^{ i < \xi, X_t>} \Phi ( X_t - y_0 ) \right) \\
= E_x \left( e^{ i < \xi, X_t>} \Phi ( X_t - y_0 )1_{\Phi ( X_t - y_0 ) > 0 ; \sup_{ 0 \le s \le \delta } | \XXX_s ( X_{\tau_1}) - X_{ \tau_1} | \geq R }  \right)\\
+ E_x \left( e^{ i < \xi, X_t>} \Phi ( X_t - y_0 )1_{ \Phi ( X_t - y_0 ) > 0 ; t- \delta = \tau_1 < t < \tau_2} \right) .
\end{multline}
The first term can be controlled, for all $ q > 0 , $ as follows.
\begin{equation}\label{eq:control11}
P_x\left( \Phi ( X_t - y_0 ) > 0 ; \sup_{ 0 \le s \le \delta } | \XXX_s ( X_{\tau_1}) - X_{ \tau_1} | \geq R \right) \le 
C(T,q,m , b, \sigma ) R^{-q} \delta^{q /2}  .
\end{equation}
Here we have used the following classical estimate: 
For all $ 0 \le s \le t \le T ,$
\begin{equation}\label{eq:ub1}
E\left( \sup_{ r :s \le r \le t } | \bar X_r^i - \bar X_s^i |^q \right) \le C(T,q ,m,  b ,  \sigma ) (t-s)^{q/2} .
\end{equation}

The above estimation in (\ref{eq:control11}) holds uniformly in $x.$ The constant $C(T,q ,m,  b ,  \sigma ) $ depends  on the supremum norms of $\bar b $ and $\bar \sigma,$ hence, by construction, on the supremum norms of $b$ and $\sigma $ on $B_{5R} (y_0 ) .$ 

The important contribution comes from the second term which can be controlled as follows.  
\begin{equation}\label{eq:control2}
E_x \left( e^{ i < \xi, X_t>} \Phi ( X_t - y_0 )1_{ \Phi ( X_t - y_0 ) > 0 ; t - \delta = \tau_1 < t < \tau_2} \right) \le 
 \sup_{ y \in B_{3R} (y_0)} | E \left( e^{ i < \xi, \XXX_\delta (y) > } \Phi ( \XXX_\delta (y)  - y_0 ) \right) | .
\end{equation}
Here, we have used the Markov property with respect to the time $ t - \delta .$ Again this control holds uniformly in $x.$ To the last term in (\ref{eq:control2}) we apply the integration by parts formula of Malliavin's calculus. We derive two times with respect to each space variable, i.e. we define the multi-index
$$ \beta = ( 1,1, 2,2, \ldots , m,m ) .$$

Then, since $ \partial_{x_k} e^{ i < \xi, x >}  = i \xi^k e^{ i < \xi , x > } , $
\begin{multline*}
| E \left( e^{ i < \xi , \XXX_\delta (y) >} \Phi ( \XXX_\delta (y)  - y_0) \right) | \le \\
\frac{ 1}{\prod_{i =1 }^m | \xi^i |^2 } | E \left( \partial_\beta e^{ i < \xi, \XXX_\delta (y)>} \Phi ( \XXX_\delta (y)- y_0 ) \right) | \le \\
\frac{ 1}{\prod_{i=1 }^m | \xi^i |^2 } E \left( | H_\beta ( \XXX_\delta (y) , \Phi ( \XXX_\delta (y) - y_0 ) ) | \right ) .
\end{multline*}
The last inequality follows from the integration by parts formula of Malliavin's calculus, and $H_\beta $ is the weight appearing in this formula, see e.g. Proposition 2.1 of De Marco 2011. We have $ | \beta | = 2 m .$ Recall that $L $ is the number of brackets needed in order to span $\bbr^m $ in $y_0,$ see condition $({\bf H3} ).$ We will show in the appendix that the following classical result holds. There exists a constant $n_L $ such that 
\begin{equation}\label{eq:classical}
\|  H_\beta ( \XXX_\delta (y) , \Phi ( \XXX_\delta (y) - y_0 )\|_p \le C (T,  p, R, m)  \delta^{ - m n_L  } .
\end{equation}
We deduce from (\ref{eq:control11}) and (\ref{eq:classical}) that, for any $ q \geq 1,$ 
$$
m_0 \hat \nu (\xi ) \le C (T, q, R, m ) \;
 \left[ R^{- q } \delta^{ q/2} + \frac{1}{\prod_{ i = 1}^m | \xi^i |^2 } \delta^{ - m n_L } \right] .
$$

The following argument is the main idea of balance that is given in De Marco 2011: We choose for a given $\xi$  a value of $\delta$ ensuring that $ R^{- q } \delta^{ q/2} + \frac{1}{\prod_{ i = 1}^m | \xi^i |^2 } \delta^{ - m n_L }$ tends to zero faster than $ \left( \prod_{i=1}^m  | \xi^i | \right)^{ - 3/2  } .$ Let $ \| \xi \| := \prod_{i=1}^m  | \xi^i |  $ and choose 
$$ \delta = t/2 \wedge 1 \wedge \| \xi \|^{-\frac{1}{2mn_L} }, \; q = 6 m  n_L  .$$
With this choice we have
\begin{equation}\label{eq:ub5} 
m_0 \hat \nu (\xi ) \le C (T, q, R, m ) \; \| \xi \|^{ - \frac32} ,
\end{equation}
and this is integrable in $\xi $ for $\| \xi \| \to \infty.$

Now we can conclude the proof of Theorem \ref{theo:main2}. Recall the definition of $\nu $ in (\ref{eq:mu}). Then for any $ y \in B_R (y_0), $ 
\begin{equation}\label{eq:ouf}
p_{0, t } ( x, y ) = \frac{m_0}{(2 \pi)^m} \int_{\bbr^m} e^{ - i < \xi , y > } \hat \nu (\xi ) d \xi =
\frac{1}{(2 \pi)^m} \int_{\bbr^m } e^{ - i <\xi , y > } E_x ( e^{ i < \xi , X_t>} \Phi ( X_t - y_0 ) ) d \xi .
\end{equation}
We cut the above integral into the integral over a finite region $I$ where $ \| \xi \| \le C $ and its complementary. 
On $I,$ we can upper bound the integrand by $ 1 $ (recall that $ \Phi \le 1_{ B_{2R} (0)} $), and on $I^c ,$ we use the above upper bound (\ref{eq:ub5}). This proves the continuity of $ p_{ 0, t} (x, y ) $ with respect to $y.$ Note that this continuity is uniform in $x,$ since the upper bounds obtained in (\ref{eq:control11}) and (\ref{eq:control2}) do not depend on the starting point $x.$

It remains to prove the lower semi-continuity of $p_{0, t} (x,y) $ in $x \in E ,$ for fixed $y \in B_R (y_0).$ The idea is to compare the diffusion $X$ to an approximation $X^n, $ which is obtained when considering $X$ before the first exit time of $K_n, $ for some fixed compact $K_n.$ It is then natural to use the flow property of $X^n$ which implies continuous dependence on the starting point. (Notice that the process $X$ itself might not satisfy the flow property.) 

For that sake, fix $n$ and let $b^n (t,x) $ and $\sigma^n (x)$ be $C^\infty_b-$extensions (in $x$) of $ b(t, \cdot _{ | K_n} )$ and $\sigma_{| K_n} .$ Let 
$X^n $ be the associated diffusion process. $X$ coincides with $X^n$ up to the first exit time $T_n.$
Hence, for $x \in K_n,$ we can write
\begin{eqnarray*}
E_x \left( e^{ i < \xi , X_t >} \Phi ( X_t - y_0) \right) &=& \lim_{n \to \infty } E_x \left( e^{i < \xi, X^n_t>} \Phi (X_t^n - y_0 ) 1_{ \{ T_n > t \} } \right) \\
& \geq & E_x \left( e^{i < \xi, X^n_t>} \Phi (X_t^n - y_0 )  1_{ \{ T_n > t \}} \right) .
\end{eqnarray*}
Here, the first equality follows from $ T_n \to \infty $ almost surely as $n \to \infty.$ The inequality follows from the fact that $X_t^n = X_t$ on $\{ T_n > t \} .$
The last expression $ E_x \left( e^{i < \xi, X^n_t>} \Phi (X_t^n - y_0 )  1_{ \{ T_n > t \}} \right)$ would depend continuously on $x$, due to the Feller property of $X^n ,$ if there would not be the presence of the indicator of $ \{ T_n > t \} .$ We have to approach the above indicator function by some continuous operation on the space of continuous functions. To be more precise, let $ \Omega = C( \bbr_+ , \bbr^m ) .$ We endow $\Omega $ with the topology of uniform convergence on compacts and write $\mathbb{P}_{0,x}^n$ for the law of $X^n $ on $ (\Omega , {\cal B} (\Omega ) ),$ starting from $x$ at time $0.$ Then we know that the family of associated probability measures $\{ \mathbb{P}_{0,x}^n, x \in \bbr^m \}$ is Feller, i.e. if $x_k \to x, $ then $\mathbb{P}_{0, x_k}^n \to \mathbb{P}_{0,x}^n $ weakly as $k \to \infty .$  

What follows is only devoted to replace the indicator of $ \{ T_n > t \} $ by some continuous functional on $ \Omega .$ Let $M_t^n = \max_{ s \le t } X_s^n $ and $m_t^n = \min_{s \le t} X_s^n $ be the (coordinate-wise) maximum and minimum processes associated to $X^n.$  Due to the structure of the compacts $K_n,$ we can construct  
$C^\infty_b-$functions $ \varphi^n , \Phi^n $ such that $ 1_{ [a_{n-1}, \infty [ } \le \varphi^n \le 1_{ [a_n,  \infty[ } $ and 
$ 1_{ ] - \infty , b_{n-1}] } \le \Phi^n \le 1_{ ] - \infty , b_n] } $ (these inequalities have to be understood coordinate-wise). Then, since $X_t$ equals $X_t^n $ up to time $T_n,$ 
$$ \{ T_{n-1} > t \} = \{ a_{n -1} \le m_t^n \le M_t^n \le b_{ n - 1} \} \subset \{ \varphi^n ( m_t^n ) > 0 , \Phi^n (M_t^n ) > 0  \} \subset \{ T_n > t \} . $$ 
So
\begin{eqnarray*}
E_x \left( e^{i < \xi, X^n_t>} \Phi (X_t^n - y_0 ) 1_{ \{ T_n > t \} } \right) 
& \geq & E_x \left( e^{i < \xi, X^n_t>} \Phi (X_t^n - y_0 )  \Phi^n ( M_t^n) \varphi^n (m_t^n ) \right) ,
\end{eqnarray*}
for any $ n .$ Write
$$ \gamma_n ( x, \xi) := \frac{1}{m_0}  E_x \left( e^{i < \xi, X^n_t>} \Phi (X_t^n - y_0 )  \Phi^n ( M_t^n) \varphi^n (m_t^n ) \right) .$$

By the Feller property of $\mathbb{P}_{0,x}^n $ and since all operations appearing in $ \gamma_n (x, \xi) $ are continuous operations on $\Omega , $ $\gamma_n ( \cdot , \xi)$ is continuous in $x,$ for any fixed $n.$ Now, instead of applying Malliavin calculus to $E_x \left( e^{ i < \xi , X_t > } \Phi ( X_t - y_0) \right) $ as we did in (\ref{eq:fourier}) above, we apply the above estimates to $ \gamma_n ( \cdot, \xi ) .$ 
Note that the upper bounds (\ref{eq:control11}), (\ref{eq:control2}) and (\ref{eq:ub5}) hold also for $m_0 \gamma_n (x, \xi).$ Moreover, they hold uniformly in $x .$ This implies, by dominated convergence, that for any $y \in B_R ( y_0),$
$$ p_{0,t}^n (x,y) := \frac{m_0}{(2\pi)^m} \int_{\bbr^m } e^{ - i < \xi , y >} \gamma_n ( x, \xi ) d \xi $$
is continuous in $x.$ Finally, we have that 
$$ p_{0, t} (x,y) = \lim_n \uparrow p_{0, t}^n (x,y) .$$
This implies the result, since the limit of a growing sequence of continuous functions is lower semi-continuous, and finishes the proof of Theorem \ref{theo:main2}. \halmos 

\subsection{Theorem \ref{theo:main2} implies Theorem \ref{theo:main}}\label{sec:4.4}
We check conditions {\bf (H1), (H2) } and {\bf (H3)} for ($\xi \! \!$ HH). Condition {\bf (H3)} is satisfied once Theorem \ref{theo:hoer} is proved. 

We now show that condition {\bf (H1)} is satisfied. By our assumptions, a strong solution $\xi_t$ of (\ref{eq:xi}) exists. Moreover, the coefficients of $V$ and $n,m, h$ are locally Lipschitz continuous. This implies the existence of a unique strong solution of ($\xi \! \!$ HH) which is a maximal solution, i.e. exists up to some explosion time. So all we have to do is to prove that the process does not explode. By assumption, $\xi_t$ does not explode. Consider now the unique solution $(V_t, n_t, m_t, h_t, \xi_t) $ of ($\xi \! \!$ HH) on $ [ 0, T_\infty [, $ where $T_\infty $ is the associated explosion time. We show first 
that $n,m$ and $h$ stay in $ (0,1) ,$ whenever they start in $(0, 1) .$ The result is a consequence of the common structure of the equations satisfied by $n,m,h$. The details are given for $n$ but the same arguments apply to $m$ and $h$. We fix $\omega $ and rewrite 
\begin{eqnarray*}
dn_t&=& (\alpha_n(V_t)(1-n_t)-\beta_n(V_t)n_t) dt \\&=& 
- ( \alpha_n + \beta_n) (V_t) n_t dt + \alpha_n ( V_t) dt \\
&=& [-a(V_t)n_t+b(V_t)]dt,
\end{eqnarray*}
where $a(v) = (\alpha_n + \beta_n) (v) $ and $ b(v) = \alpha_n ( v ).$
Given the fixed trajectory $V_t $ on $ [ 0, T_\infty [ , $ the variation of constants method yields the following representation of $n_t :$ 
\begin{equation}\label{eq:exact}
n_t=n_0{\rm e}^{-\int_0^t a(V_s)ds}+\int_0^t b(V_u){\rm e}^{-\int_u^t a(V_r)dr}du, \; t < T_\infty .
\end{equation}
Notice that the above equation does not provide an explicit formula for $n_t,$ since $V  $ depends on $n .$  

We rewrite $\int_0^t b(V_u){\rm e}^{-\int_u^t a(V_r)dr}du=\int_0^t \frac{b(V_u)}{a(V_u)}a(V_u){\rm e}^{-\int_u^t a(V_r)dr}du$. By definition, $a(v)$ is positive and  $\frac{b(v)}{a(v)}\in [0,1).$ Hence 
\begin{equation}
0 <  n_t\leq n_0{\rm e}^{-\int_0^t a(V_s)ds}+\int_0^t a(V_u){\rm e}^{-\int_u^t a(V_r)dr}du.
\end{equation}
In the above formula, the strict positivity of $n_t$ follows from the fact that $ \int_0^t a(V_s ) ds + \int_0^t b(V_s) ds < \infty, $ since $t < T_{\infty } .$ By integration by parts it follows that
\begin{eqnarray}\label{eq:tobecited}
0 <  n_t&\leq& (n_0+{\rm e}^{\int_0^t a(V_r)dr}-1){\rm e}^{-\int_0^t a(V_s)ds} \nonumber \\
&=& 1+(n_0-1){\rm e}^{-\int_0^t a(V_s)ds}.
\end{eqnarray}
Therefore if $n_0\in (0,1)$, then for all $t < T_\infty $, $n_t\in (0,1)$. The same kind of arguments apply to $m$ and to $h.$ 

As a consequence, we deduce immediately from ($\xi \! \!$ HH) that for suitable constants $C_1 $ and $C_2,$
$$ | V_t | \le | \xi_t | + C_1 \int_0^t | V_s| ds + C_2 t , \; t  < T_\infty .$$ 
This implies, using Gronwall's inequality and non explosion of $\xi_t,$ that $V_t$ does not explode neither. 
Hence $T_\infty = \infty $ almost surely and the above estimates hold on $ [ 0, \infty [ .$ 

Now, let $C_n \subset C_{n+1} \subset \bbr$ be a growing sequence of compact intervals such that $\bigcup C_n = U .$ Put
$$ K_n = [ - n, n  ] \times [  \frac1n, 1 - \frac1n ]^3 \times C_n ,$$
then $ T_n = \inf \{ t : X_t \in K_n^c\} \to \infty $ as $ n \to \infty .$ Moreover,
clearly {\bf (H2)} is satisfied on $K_n.$ Therefore, all conditions needed in order to apply Theorem \ref{theo:main2} are satisfied, and thus Theorem \ref{theo:main} follows.

\section{Proof of Theorem 1.}\label{section:th1}
Let $X_t = (V_t, n_t, m_t, h_t, \xi_t )$ be the $5-$dimensional diffusion process of ($\xi \! \!$ HH). Write 
$$ b (t,x)= \left( 
\begin{array}{c}
b^1 (t,x) \\
\vdots \\
b^5 (t,x) 
\end{array}
\right) \in \bbr^5  \mbox{ and } \sigma (x)= \left( 
\begin{array}{c}
\sigma^1 (x) \\
\vdots \\
\sigma^5 (x) 
\end{array}
\right) \in \bbr^5,  x = ( v, n, m , h , \zeta ) \in E_5 ,$$ for its drift function and its diffusion coefficient. Here, 
$$ b^1 (t,x) = (S(t) - \zeta ) \tau  - F( v , n , m, h) ,\; b^5 (t,x) = (S(t) - \zeta ) \tau, $$
$$ b^2 (t, x) = G_n ( v, n), \; b^3 ( t, x) = G_m ( v, m), \; b^4 (t, x) = G_h ( v,h) .$$
Moreover, writing 
\begin{equation}\label{eq:a}
d(\zeta ) := \gamma \sqrt{\tau} q(\zeta) ,
\end{equation} 
we have
$ \sigma^1 (x) = \sigma^5 (x) = d(\zeta)   $
and 
$ \sigma^2 (x) = \sigma^3 (x) = \sigma^4 (x) = 0 .$ Hence ($\xi   $HH) can be rewritten as five dimensional diffusion equation
$$ dX_t = b (t, X_t) dt + \sigma (X_t) d W_t.$$
As before, we rewrite this equation in the Stratonovitch sense and introduce $\tilde b (t,x),$ where $ \tilde b^i (t,x) = b^i ( t, x) , $ for $i = 2 ,3, 4,$ and 
\begin{equation}\label{eq:strato}
 \tilde b^i (t,x) = b^i (t,x) - \frac12 d' (\zeta) d(\zeta)  , i = 1 , 5 .
\end{equation}
Since this drift is time dependent, the associated directional derivative is 
$$ \frac{\partial }{\partial t} + \tilde b^1 (t,x) \frac{\partial}{\partial v} + \ldots + \tilde b^5 (t,x) \frac{\partial }{\partial \zeta} .$$
We start by calculating the Lie-bracket of $ \sigma $ and $\tilde b .$ In order to simplify notation, we identify the vector $x= (v,n,m,h, \zeta) $ with $x= (x^1 , \ldots , x^5 ).$ Then
\begin{eqnarray*}
 [ \tilde b , \sigma ]^i & = &\frac{\partial \sigma^i }{\partial t} + \sum_{ j= 1 }^5 \left( \tilde b^j \frac{\partial \sigma^i }{\partial x^j } - \sigma^j \frac{\partial \tilde b^i }{\partial x^j } \right)  =\sum_{ j= 1 }^5 \left( \tilde b^j \frac{\partial \sigma^i }{\partial x^j } - \sigma^j \frac{\partial \tilde b^i }{\partial x^j } \right)\\
&=&  d' (\zeta) \tilde b^5 (t,x) (\delta_{i1} + \delta_{i5}) - d(\zeta) \left( \frac{\partial \tilde b^i }{\partial v} + \frac{\partial \tilde b^i }{\partial \zeta} \right)  .
\end{eqnarray*}
As a consequence, we get 
\begin{eqnarray*}
 [ \tilde b , \sigma ]&=& - d(\zeta) \left( 
\begin{array}{c}
-\tau - \frac12  [ (d'(\zeta))^2 + d(\zeta) d'' (\zeta) ]- \partial_v F(v,n,m,h) \\
 g_n(v,n)\\
 g_m (v,m)\\
  g_h ( v,h)\\
  -\tau - \frac12  [ (d'(\zeta))^2 + d(\zeta) d'' (\zeta) ]
\end{array}
\right) 
+  \left( 
\begin{array}{c}
d' (\zeta) \tilde b^5 (t,x)\\
0\\
0\\
0\\
d' (\zeta) \tilde b^5 (t,x)
\end{array}
\right) \\
&=& d(\zeta) \left( 
\begin{array}{c}
\partial_v F(v,n,m,h) \\
- g_n(v,n)\\
- g_m (v,m)\\
 - g_h ( v,h)\\
0
\end{array}
\right) 
+ A_2 (t, \zeta) \left( 
\begin{array}{c}
1 \\
0\\
0\\
0\\
1 
\end{array}
\right)
,
\end{eqnarray*}
where 
$$ A_2 (t, \zeta ) :=  d' (\zeta) \tilde b^5 (t,x) + d(\zeta) [ \tau + \frac12  [ (d'(\zeta))^2 + d(\zeta) d'' (\zeta) ] ]  .$$
Write $ V_{ 2} =  [ \tilde b , \sigma ] .$ We are now going to evaluate $ V_{3} = [ \sigma , V_{2}] .$ We have
\begin{eqnarray*}
 [\sigma , V_{2} ]^i &=& \sum_{j = 1}^5 \left( \sigma^j \frac{ \partial V_{2}^i }{\partial x^j } - V_{2}^j \frac{\partial  \sigma^i }{\partial x^j } \right) = a (\zeta) \left( \frac{\partial V^i_{2}}{\partial v} + \frac{\partial V^i_{2}}{\partial \zeta}  \right) - V_{2}^5 \frac{\partial \sigma^i }{\partial \zeta} .
\end{eqnarray*}
Therefore, 
$$
V_3 (t,x) = d^2 (\zeta) \left( 
\begin{array}{c}
\partial^2_v F (v,n,m,h) \\
- g_n' (v,n)\\
- g_m' (v,m)\\
-g_h' (v,h) \\
0 
\end{array}
\right) + 
 d(\zeta) d' (\zeta) 
\left( 
\begin{array}{c}
\partial_v F(v,n,m,h) \\
- g_n(v,n)\\
- g_m (v,m)\\
 - g_h ( v,h)\\
0
\end{array}
\right) 
+A_3 ( t, \zeta)  
\left( 
\begin{array}{c}
1\\
0\\
0\\
0\\
1
\end{array}
\right) ,
$$
where
$$ A_3 (t,\zeta ) :=  d(\zeta)  \partial_\zeta A_2 (t,\zeta) - A_2 (t, \zeta)  d' (\zeta ) .$$
Putting $V_4 = [ \sigma, V_3 ] ,$ we obtain analogously that 
\begin{multline*}
V_4 = d^3 (\zeta) \left( 
\begin{array}{c}
\partial^3_v F(v,n,m,h) \\
- g_n'' (v,n)\\
-g_m'' (v,m)\\
- g_h'' (v,h) \\
0 
\end{array}
\right) 
+ \left( 3 d^2 (\zeta) d' (\zeta) \right) 
\left( 
\begin{array}{c}
\partial^2_v F(v,n,m,h) \\
- g_n' (v,n)\\
-g_m' (v,m)\\
- g_h' (v,h) \\
0 
\end{array}
\right) 
\\
+ \left(  d' ( \zeta)^2 + d (\zeta ) d'' (\zeta)  \right) 
\left( 
\begin{array}{c}
\partial_v F(v,n,m,h) \\
- g_n (v,n)\\
-g_m (v,m)\\
- g_h (v,h) \\
0 
\end{array}
\right)
+  A_4 ( t, \zeta) 
\left( 
\begin{array}{c}
1\\
0\\
0\\
0\\
1
\end{array}
\right) ,
\end{multline*}
where 
$$ A_4 ( t, \zeta) =  d(\zeta) \partial_\zeta A_3 (t,\zeta) - A_3 (t, \zeta ) d' (\zeta)  .$$
Finally, for $ V_5 = [ \sigma , V_4] ,$ we obtain similarly a representation
$$V_5 =  d^4 (\zeta) 
\left( 
\begin{array}{c}
\partial^4_v F(v,n,m,h) \\
- g_n''' (v,n)\\
-g_m''' (v,m)\\
- g_h''' (v,h) \\
0 
\end{array}
\right)
+ a_1 (\zeta) B_1 + a_2 (\zeta) B_2 + a_3 (\zeta) B_3 + A_5 (t,\zeta)\left( 
\begin{array}{c}
1\\
0\\
0\\
0\\
1
\end{array}
\right) ,$$
where 
$ a_1, a_2, a_3 $ are functions defined in terms of derivatives of $d$ and where $ B_1, B_2, B_3 \in Span \{ V_2, V_3, V_4 \} .$   

Now we are able to conclude our proof. Since by definition of $F$ in (\ref{eq:F}), $ \partial_v F (v,n,m,h) \neq 0 $ for all $ (v,n,m,h) \in E_4$ and $ \partial^k_v F(v,n,m,h) \equiv 0 $ for all $ k \geq 2 ,$  we have for all fixed $x \in E_5,$  
$$ Span \{ \sigma , V_2 , V_3, V_4, V_5 \} = Span \{ e_1, e_5, g_1, g_2, g_3 \}, $$
where 
$$ g_i = \left( 
\begin{array}{c}
0\\
\partial^i_v g_n ( v,n)\\
\partial^i_v g_m (v,m)\\
\partial^i_v g_h (v,h) \\
0 
\end{array}
\right) .$$
By the definition of $D ( v,n,m,h) $ in (\ref{eq:det}), these five vectors $\sigma , V_2, V_3 , V_4 $ and $V_5$ span $\bbr^5 $ for all $x \in e_5$ such that $D (v, n,m,h) \neq 0 .$ This concludes the proof. \halmos

\section{Proof of Theorem 3.}\label{sec:th3}

With all notations of Theorem 3, we consider the system ($\xi$HH) driven by $S$ of Section \ref{sec:2.5}, 
\begin{equation}
X_s = x + \int_0^s \sigma ( X_u)  dW_u + \int_0^s  b( u, X_u) du , \quad s\le t .
\end{equation}
We write 
$ \Omega  = C ( [ 0, \infty [ , \bbr^5 )  $
and endow $\Omega $ with its canonical filtration $ ( {\cal F}_t)_{t \geq 0 } .$ Let $ \mathbb{P}_{0,x}$ be the law of $(X_{  u }, u \geq 0 ), $ starting from $x.$ In order to find lower bounds for quantities of the form $\mathbb{P}_{0,  x } (  B ) $ for measurable $B \in {\cal F}_t ,$ we will use control arguments and the support theorem for diffusions. We need first to localize the system.  
Let $C_n \subset C_{n+1 } \subset U $ be a sequence of compact intervals such that $ \bigcup C_n = U .$ Put $K_n = [ - n , n ] \times [ \frac1n , 1 - \frac1n ]^3 \times C_n \subset E_5$ and let $T_n = \inf \{ t : X_t \in K_n^c \} $ be the exit time of $K_n.$ For a fixed $n,$ let $ b^n (t, x) $ and $\sigma^n (x)$ be $C_b^\infty -$extensions in $x$ of $ b(t, \cdot_{| K_n }) $ and $\sigma_{| K_n} .$  
Let $X^n$ be the associated diffusion process. Then for any starting point $x \in K_n ,$ we write $\mathbb{P}_{0, x}^n $ for the law of $(X^n_{  u } , u \geq0  )$ on $\Omega . $ Fix a time $t > 0 .$ Then for any measurable $B \in {\cal F}_t ,$ 
\begin{equation}\label{eq:tobelb}
\mathbb{P}_{0,  x } (  B ) \geq \mathbb{P}_{ 0, x } ( \{ f  \in B ; T_n > t \}  ) 
 = \mathbb{P}^n_{ 0, x } (  \{ f \in B ; T_n >  t \} )  .
\end{equation}
It suffices to show that this last expression is strictly positive, for suitable choices of $x \in K_n$ and $B.$ 
For this sake, as already mentioned,  we will use the support theorem for diffusions, see Stroock and Varadhan 1972.  Let $ {\cal H} = \{ h : [ 0, t ] \to \bbr : h(s) = \int_0^s \dot h (u) du , \forall s \le t , \int_0^t \dot h^2 (u) du < \infty \} $ be the Cameron-Martin space. 
Given $h \in {\cal H} ,$ consider $ X(h)$ the solution of the differential equation 
\begin{equation}\label{eq:control1}
X(h)_s = x + \int_0^s \sigma^n ( X(h)_u) \dot h (u) du + \int_0^s \tilde b^n ( u, X(h)_u) du , \quad s\le t ,
\end{equation}
where in accordance with the notation used previously in the paper, $X(h)$ is $5-$dimensional of the form $X(h)=(X(h)^1, X(h)^2, X(h)^3, X(h)^4, S(h))$. In the above formula (\ref{eq:control1}), $\tilde b^n$ is the drift vector of $X^n,$ written in Stratonovitch form. 

As a consequence of the support theorem for diffusions (see e.g. Theorem 3.5 of Millet and Sanz-Sol\'e 1994 or Theorem 4 of Ben Arous, Gradinaru and Ledoux 1994), the support of the law $\mathbb{P}^n_{0, x}$ in restriction to $ {\cal F}_t$ is the closure 
of the set $ \{ X(h) : h \in {\cal H} \} $ with respect to the uniform norm on $ [ 0, t] .$  

In order to find lower bounds for (\ref{eq:tobelb}) we have to construct solutions $X(h)$ of (\ref{eq:control1}) which stay in $K_n$ during $ [ 0, t ] . $ But on $K_n,$ both processes $X^n$ and $X$ have the same coefficients. Hence, in restriction to $K_n,$ the above control problem (\ref{eq:control1}) is equivalent to the following, where we recall that $d(.)$ has been defined in (\ref{eq:a}). 
$$
\left\{\begin{array}{l}
\frac{d}{ds }  X(h)^1_s   \;=\;  \frac{d}{ds} S(h)_s \;- F( X(h)^1_s, X(h)^2_s, X(h)^3_s, X(h)^4_s) \\
\frac{d}{ds }  X(h)^2_s \;=\;  \, \al_n(X(h)^1_s)\,(1-X(h)^2_s)  \;-\; \beta_n(X(h)^1_s)\, X(h)^2_s  \,  \\
\frac{d}{ds }  X(h)^3_s\;=\;  \, \al_m(X(h)^1_s)\,(1-X(h)^3_s)  \;-\; \beta_m(X(h)^1_s)\, X(h)^3_s  \,\\
\frac{d}{ds }  X(h)^4_s \;=\; \, \al_h(X(h)^1_s)\,(1-X(h)^4_s)  \;-\; \beta_h(X(h)^1_s)\, X(h)^4_s  \,  \\
\frac{d}{ds }  S(h)_s\;=\; (\, S(s)-S(h)_s\,)\, \tau - \frac{1 }{2} d' ( S(h)_s) d( S(h)_s)  \;+\;  \gamma q (S(h)_s)\, \sqrt{\tau}\;  \dot h (s)  .
\end{array}\right. 
\leqno{\rm (HHcontrolled)}
$$

In order to find simple solutions of the above system, we consider the specific starting point $x  $ which is prescribed in Theorem \ref{theo:pos}.  Since $\zeta  \in U $ and $ (v,n,m,h) \in E_4,$ there exists $n$ such that $ x = ( v,n,m,h, \zeta  ) \in K_n .$ 

We will use $x$ as starting point and construct solutions of ${\rm (HHcontrolled)}$ such that 
\begin{equation}\label{eq:sol}
\frac{d}{ds }  S(h)_s = \tilde S (s)  , \quad \mbox{ for all $s \le t.$ }
\end{equation}
Equation (\ref{eq:sol}) implies that $ S(h)_s = S(h)_0+ \int_0^s \tilde S(u) du = \zeta + \int_0^s \tilde S(u) du = : \tilde I_s, $ for all $s \le t.$  
Hence, if we define
$$ \dot h (s) := \frac{ \tilde S(s) + (\tilde I_s - S(s) ) \tau + \frac{1 }{2} d' (  \tilde I_s) d( \tilde I_s)  }{ \gamma q (\tilde I_s ) \sqrt{\tau}} ,  \quad  s\le t ,$$
then the right hand side of the last line of ${\rm (HHcontrolled)}$ equals indeed $\tilde S (s) .$ 

Notice that $\dot h $ is well-defined since $ q(\tilde I_s) > 0 $ for all $s \le t.$ Moreover, the signals $S$ and $\tilde S$ being $T-$periodic, clearly $\dot h \in L^2 ( [ 0, t ] ) ,$ hence $ h \in {\cal H} .$ With this choice of $h,$ the first four lines of ${\rm (HHcontrolled)}$ reduce to the deterministic system (HH)  with input signal $s \to \tilde S(s) .$ Write $\mathbb{Y}$ for the associated deterministic solution starting from $ (v, n,m,h) $ at time $0$ and $\mathbb {X}_s = ( \mathbb{Y}_s, \tilde I_s ), s \le t .$ Then for $n$ sufficiently large, $ \mathbb{X}_s \in K_n $ for all $ s \le t .$

By the support theorem, for every $\varepsilon > 0, $ putting $B^\infty_\varepsilon ( \mathbb{X} ) = \{ f \in  \Omega  : \sup_{ s \le t } | f(s) - \mathbb{X}_s | <  \varepsilon \} ,$ we have that 
$$ \mathbb{P}^n_{0, x} (B_\varepsilon^\infty (\mathbb{X}) ) > 0 .$$
 
Now choosing $ \varepsilon $ such that $ B^\infty_\varepsilon ( \mathbb{X} )  \subset \{ f \in \Omega : T_n (f) > t \} $ and putting 
$ B=  B^\infty_\varepsilon ( \mathbb{X} )  ,$ we obtain the desired first result of Theorem \ref{theo:pos}. 
Finally, by the Feller property of $ \mathbb{P}^n_{0, x } ,$ for fixed $\varepsilon ,$ we can extend the above property to a small ball around $x .$ This shows the second assertion of Theorem \ref{theo:pos}. \halmos

\section{Appendix : Some elements of Malliavin calculus}
In this appendix we give the basic arguments from Malliavin calculus that allow to show that the important estimate (\ref{eq:classical}) holds true in the non time homogeneous case as well as in the time homogeneous case. For the basic concepts of Malliavin calculus, we refer the reader to the classical reference Nualart 1995. 

Throughout this section, $\bar X  $ denotes the unique strong solution of the SED (recall also (\ref{eq:processgood}))
\begin{equation}\label{eq:processgoodbis}
\bar X_t^i = x^i + \int_0^t \bar b^i ( s, \bar X_s) ds + \int_0^t \bar \sigma^i (\bar X_s) d W_s , t \le T , 1 \le i \le m ,
\end{equation}
where $x \in \bbr^m ,$ $ \bar \sigma^i \in C^\infty_b ( \bbr^m ) ,  $  and where $\bar b^i (t,x) $ and all partial derivatives $ \partial^\alpha_x \partial^\beta_t \bar b^i (t, \cdot ) $ are bounded uniformly in $t \in [0, T ] .$ 
Let 
$$ \tilde{\bar b}^i (t,x) = \bar b^i (t,x) - \frac12 \sum_{k=1}^m \bar\sigma^k (x) \frac{ \partial \bar \sigma^i }{\partial x_k } (x) ,\;  1 \le i \le m ,$$
and $\bar V_0  = \frac{\partial}{\partial t } + \tilde{\bar b}  $ be the associated time-space directional derivative. We use analogous notation to section \ref{sec:hoerm} und put $ \bar V_1 = \bar \sigma .$
The local H\"ormander condition for a point $x$ and a given number of brackets $L$ is 
$$  \bar {\cal V}_L (x) > 0 , \leqno{ {\rm \bf (H_L )}(x)} $$
where $ \bar {\cal V}_L (x) $ is defined analogously to (\ref{eq:lie1}). 

The main ingredient for the control of the weight in Malliavin's integration by parts formula as in formula (\ref{eq:classical}) is to obtain estimates of Malliavin's covariance matrix. We check that all results obtained in Kusuoka- Stroock 1985 are still valid in our framework. Let 
$$ (Y_t)_{i,j} = \frac{ \partial \bar X_t^i }{\partial x_j} , 1 \le i, j \le m .$$
Then $Y$ satisfies the following linear equation having bounded coefficients (bounded with respect to time and space)
$$ Y_t = I_m + \int_0^t  \partial \bar b ( s, \bar X_s) Y_s  ds + \int _0^t \partial \bar \sigma ( \bar X_s) Y_s d W_s .$$
Here $I_m$ is the $m\times m-$unity matrix and $ \partial \bar b $ and $ \partial \bar \sigma $ are the $m\times m-$matrices having components $ ( \partial \bar b)_{i,j} (t, x) = \frac{\partial \bar b^i }{\partial x_j} (t, x ) $ and  $ ( \partial \bar \sigma)_{i,j} ( x) = \frac{\partial \bar \sigma^i }{\partial x_j} ( x )  .$
By means of It\^o's formula, one shows that $Y_t$ is invertible. The inverse $Z_t$ still satisfies a linear equation with coefficients bounded in $t$ and in $x$ given by 
\begin{equation}\label{eq:z}
Z_t = I_m - \int_0^t \partial \tilde{\bar b} (s,\bar X_s)  Z_s ds - \int_0^t  \partial \bar \sigma ( \bar X_s) Z_s \circ dW_s   ,
\end{equation}
where $ \circ d W_s $ denotes the Stratonovitch integral. In this framework, the following estimates are classical (see e.g.  Kusuoka and Stroock 1985 or De Marco 2011, Prop. 2.2 and Lemma 2.1). 
For all $ 0 \le s \le t \le T ,$ for all $ p \geq 1,$ 
\begin{equation}\label{eq:ub1}
E\left( \sup_{ r :s \le r \le t } | \bar X_r^i - \bar X_s^i |^p \right) \le C(T,p,m, \bar b , \bar \sigma ) (t-s)^{p/2} ,
\end{equation}
\begin{equation}\label{eq:ub2}
\sup_{ s \le t } E ( |(Z_s)_{i,j}|^p) \le C(T,p,m, \bar b , \bar \sigma ) ,\; 1 \le i, j \le m ,
\end{equation} 
\begin{equation}\label{eq:ub3}
\sup_{ r_1 , \ldots , r_k \le t } E \left(  | D_{r_1, \ldots , r_k} \bar X_t^i |^p  \right) \le C(T,p,m,k, \bar b , \bar \sigma )  \left( t^{1/2}  +  1\right)^{ (k+1)^2 p} ,
\end{equation}
where the constants depend only on the bounds of the derivatives with respect to space of $\bar b$ and $\bar \sigma .$ Notice that the above estimates are not sharp, and much better estimates can be obtained, see for instance in De Marco 2011. However, for our purpose, the above estimates are completely sufficient.

As indicated before, the main issue in order to prove (\ref{eq:classical}) is to obtain estimates on the Malliavin covariance matrix. 
So let $ (\sigma_{ \bar X_t})_{i,j}  = < D \bar X_t^i , D \bar X_t^j>_{L^2 [0, t ]} , 1 \le i, j \le m .$
Then it is well known, see for example formula (240), page 110 of Nualart 1995, that 
$$ \sigma_{ \bar X_t} = Y_t \left( \int_0^t Z_s \bar \sigma ( \bar X_s) \bar \sigma^* ( \bar X_s) Z_s^* ds \right) Y_t^* .$$
In order to evaluate the inner integral, one has to control expressions of the type 
$ Z_s V(s, \XXX_s) ,$ where $ V (t,x) $ is a smooth function of $t$ and $x.$ Using partial integration it is easy to see that 
\begin{multline}\label{eq:important2}
Z_t V (t, \XXX_t ) = V(0, x) + \int_0^t Z_s [ \bar \sigma , V ] ( s, \XXX_s) \circ d W_s 
+ \int_0^t Z_s  [ \frac{\partial}{\partial t} + \tilde{ \bar b} , V] (s, \XXX_s) ds \\
= V(0, x) + \int_0^t Z_s [ \bar \sigma , V ] ( s, \XXX_s) \circ d W_s 
+ \int_0^t Z_s  [ \bar V_0 , V] (s, \XXX_s) ds 
\end{multline}
(see formula (2.10) of Kusuoka and Stroock 1985), 
where we recall that 
$$ [ \frac{\partial}{ \partial t } + \tilde{\bar b} , V ]^i =  \frac{ \partial V^i }{\partial t } +\sum_{ j=1}^m \left( \tilde{\bar b}^j \frac{\partial V^i }{\partial x^j } - V^j \frac{\partial \tilde{\bar b}^i  }{\partial x_j } \right).$$
Iterating (\ref{eq:important2}) we obtain completely analogously to Theorem 2.12 of Kusuoka and Stroock 1985, for any $ L \geq 1 , $
\begin{equation}\label{eq:important1}
Z_s \bar \sigma (\XXX_s) = \sum_{ \alpha : \| \alpha \| \le L- 1 } W^{(\alpha)} (s) (\bar V_1)_{(\alpha)} (0, x ) + R_L (s, x, \bar V_1 ),
\end{equation}
where $R_L$ is a remainder term and where $ W^{(\alpha)} $ is a multiple Wiener integral. Here, $\bar V_1 = \bar \sigma $ and the $ (\bar V_1)_{(\alpha)} (0, x )$ are the successive Lie brackets. The most important feature in the above development (\ref{eq:important1}) is that the behavior of the remainder term depends only on the supremum norms of derivatives with respect to time and space of $\bar b $ and with respect to space of $\bar \sigma .$ Then, following Kusuoka and Stroock 1985, we obtain their 

\begin{cor}[Corollary 3.25 of Kusuoka and Stroock 1985]
For any  $p \geq 1 $ and $ t \le 1 , $ for any $ L \geq 1, $ for any $x$ such that (${\bf H_L})(x)$ is satisfied,
\begin{equation}\label{eq:important3}
E_x \left( | det \sigma_{\XXX_t}|^{-p} \right)^{1/p} \le C(p,m, L) \frac{1}{(\bar {\cal V}_L (x)^{ 1 + \frac{2}{L}} t )^{mL}} .
\end{equation}
\end{cor}

Once this control (\ref{eq:important3}) is established, the upper bound (\ref{eq:classical}) follows according to a well-known scheme. We refer the reader for instance to De Marco 2011, proof of Theorem 2.3., for a detailed presentation of the arguments.

\end{document}